\RequirePackage{ifthen}
\newboolean{MPA}
\setboolean{MPA}{false}

\ifthenelse {\boolean{MPA}}
{
\documentclass{svjour3}
\smartqed
\usepackage[margin=1.3in]{geometry}
} {
\documentclass[11pt]{article}
\usepackage[margin=.9in]{geometry}
}

\usepackage{graphicx}
\usepackage{amssymb}
\usepackage{amsmath}
\usepackage{verbatim,booktabs}
\usepackage{subfigure,epsfig,url,psfrag}
\usepackage{float}
\usepackage{mathrsfs}

\usepackage[usenames,dvipsnames]{pstricks}
\usepackage{epsfig}

\usepackage{cite} 

\usepackage{color}
\usepackage{enumerate}
\usepackage{enumitem}

\usepackage{todonotes}

\usepackage{hyperref}
\usepackage[capitalize,noabbrev]{cleveref}
\crefname{problem}{Problem}{Problems}
\crefname{claim}{Claim}{Claims}

\newcommand{\pare}[1]{\left(#1\right)}
\newcommand{\bra}[1]{\left\{#1\right\}}

\newcommand{\ie}{i.e., }

\newcommand{\E}{\mathcal E}

\newcommand{\I}{\mathcal I}

\def\R{{\mathbb R}}
\def\proj{{\rm proj}}

\def\01{\ensuremath{0\mathord{-}1}}
\allowdisplaybreaks

\DeclareMathOperator{\conv}{conv}
\DeclareMathOperator{\gap}{gap}
\DeclareMathOperator{\mh}{mh}
\DeclareMathOperator{\poly}{poly}
\DeclareMathOperator{\ext}{ext}
\DeclareMathOperator{\BPO}{BPO}

\DeclareMathOperator{\MS}{\mathcal S_m}
\DeclareMathOperator{\PBS}{\mathcal S_{pB}}
\DeclareMathOperator{\MP}{\mathcal P_m}
\DeclareMathOperator{\PBP}{\mathcal P_{pB}}

\ifthenelse {\boolean{MPA}}
{

\newenvironment{prf}[1][]
{\begin{proof}}
{\qed \end{proof}}

\newenvironment{prfc}[1][]
{\begin{proof}[#1]}
{\qed \end{proof}}

\newenvironment{prfh}[1][]
{\begin{proof}}
{\end{proof}}

\newcounter{claim} 
\renewenvironment{claim}[1][]
{\refstepcounter{claim} \begin{trivlist} \item[] {\bf Claim~\theclaim}\space#1 \itshape}
{\end{trivlist}}

\newenvironment{cpf}
{\begin{trivlist} \item[] {\em Proof of claim }}
{$\hfill\diamond$ \end{trivlist}}

\journalname{Mathematical Programming A}

\newtheorem{observation}{Observation}

} {

\usepackage{amsthm}
\newtheorem{theorem}{Theorem}
\newtheorem{definition}{Definition}

\newtheorem{proposition}{Proposition}
\newtheorem{example}{Example}
\newtheorem{observation}{Observation}
\newtheorem{claim}{Claim}
\newtheorem{remark}{Remark}

\newenvironment{prf}[1][]
{\begin{proof}}
{\end{proof}}

\newenvironment{cpf}
{\begin{trivlist} \item[] {\em Proof of claim. }}
{$\hfill\diamond$ \end{trivlist}}

}

\usepackage{thmtools}
\usepackage{thm-restate}

 \allowdisplaybreaks

\title{The pseudo-Boolean polytope and polynomial-size extended formulations for  binary polynomial optimization}

\ifthenelse {\boolean{MPA}}
{
\titlerunning{The pseudo-Boolean polytope and binary polynomial optimization}

\author{Alberto Del Pia \and Aida Khajavirad}

\institute{Alberto Del Pia \at
              Department of Industrial and Systems Engineering \& Wisconsin Institute for Discovery, 
              University of Wisconsin-Madison.
              E-mail: {\tt delpia@wisc.edu}.
           \and
           Aida Khajavirad \at
              Department of Industrial and Systems Engineering,
              Lehigh University.
              E-mail: {\tt aida@lehigh.edu}.
}
}
{
\author{Alberto Del Pia
\thanks{Department of Industrial and Systems Engineering \& Wisconsin Institute for Discovery,
             University of Wisconsin-Madison.
             E-mail: {\tt delpia@wisc.edu}.
             }
\and
Aida Khajavirad
\thanks{Department of Industrial and Systems Engineering,
             Lehigh University.
             E-mail: {\tt aida@lehigh.edu}.
             }
}
}

\date{July 1, 2024}

\begin{document}
\maketitle

\begin{abstract}
With the goal of obtaining strong relaxations for binary polynomial optimization problems, we introduce the pseudo-Boolean polytope defined as the set of binary points $z \in \{0,1\}^{V \cup S}$ satisfying a collection of equalities of the form 
$z_s = \prod_{v \in s} \sigma_s(z_v)$, for all $s \in S$, where $\sigma_s(z_v) \in \{z_v, 1-z_v\}$, and where $S$ is a multiset of subsets of $V$. 
By representing the pseudo-Boolean polytope via a signed hypergraph, we obtain sufficient conditions under which this polytope has a polynomial-size extended formulation. Our new framework unifies and extends all prior results on the existence of polynomial-size extended formulations for the convex hull of the feasible region of binary polynomial optimization problems of degree at least three.
\ifthenelse {\boolean{MPA}}
{
} {}
\end{abstract}

\medskip
\ifthenelse {\boolean{MPA}}
{}{
} 


\section{Introduction}
\label{sec: intro}
We consider the problem of maximizing a multivariate polynomial function over the set of binary points, henceforth referred to as~\emph{binary polynomial optimization}. 
This problem class has numerous applications across science and engineering~\cite{boros02}, and is NP-hard in general.
Based on the encoding of the polynomial objective function, we obtain two popular optimization problems, which we refer to as ``multilinear optimization'' and ``pseudo-Boolean optimization.''
We first introduce these two optimization problems in \cref{sec intro multi,sec intro Pseudo-Boolean}, and then compare and contrast them in~\cref{sec vs}.
As binary quadratic optimization has been extensively studied in the literature (see for example~\cite{barMah86,Pad89,borCraHam90}), in this paper, we focus on binary polynomial optimization problems of degree at least three.

\subsection{Multilinear optimization}
\label{sec intro multi}

A \emph{hypergraph} $G$ is a pair $(V,E)$, where $V$ is a finite set of nodes and $E$ is a set of \emph{edges}, which are subsets of $V$ of cardinality at least two.
Following the approach introduced in~\cite{dPKha17MOR}, with any hypergraph $G=(V,E)$, and cost vector $c \in \R^{V \cup E}$, we associate the following \emph{multilinear optimization problem}:
\begin{align}
\label[problem]{prob MO}
\tag{$\BPO_\mathrm{m}$}
\begin{split}
\max & \qquad \sum_{v\in V} {c_v z_v} + \sum_{e\in E} {c_e \prod_{v\in e} {z_v}} \\
{\rm s.t.} & \qquad z_v \in \{0,1\} \qquad \forall v \in V.
\end{split}
\end{align}
It can be checked that any binary polynomial optimization problem has a unique representation of the form~\cref{prob MO}.  Next, we linearize the objective function of~\cref{prob MO} by introducing a new variable $z_e$ for each product term $\prod_{v \in e} z_v$ to obtain an equivalent reformulation of this problem in a lifted space:
\begin{align}
\label[problem]{prob lMO}
\tag{$\ell\BPO_\mathrm{m}$}
\begin{split}
\max & \qquad \sum_{v\in V} {c_v z_v} + \sum_{e\in E} {c_e z_e} \\
{\rm s.t.} & \qquad z_e = \prod_{v\in e} {z_v} \qquad \forall e \in E\\
& \qquad z_v \in \{0,1\} \qquad \forall v \in V.
\end{split}
\end{align}
With the objective of understanding the facial structure of the convex hull of the feasible region of~\cref{prob lMO}, Del~Pia and Khajavirad~\cite{dPKha17MOR} introduced the~\emph{multilinear set}, defined as
\begin{equation*} \label{eq: SG}
\MS(G):= \Big\{ z \in \{0,1\}^{V \cup E} : z_e = \prod_{v \in e} {z_{v}}, \; \forall e \in E \Big\},
\end{equation*}
and its convex hull, which is called the~\emph{multilinear polytope} and is denoted by $\MP(G)$. 

In~\cite{dPKha18SIOPT,dPKha23MPA}, the authors show that the complexity of the facial structure of $\MP(G)$ is closely related to the ``acyclicity degree'' of $G$.
The most well-known types of acyclic hypergraphs, in increasing order of generality, are Berge-acyclic, $\gamma$-acyclic, $\beta$-acyclic, and $\alpha$-acyclic hypergraphs~\cite{fagin83,BeFaMaYa83,Dur12,bra14}.
We next briefly review the existing results on the facial structure of the multilinear polytope of acyclic hypergraphs. 
Recall that the \emph{rank} of a hypergraph $G$, denoted by $r$, is the maximum cardinality of any edge in $E$. 
In~\cite{dPKha18SIOPT,BucCraRod16}, the authors 
prove that $\MP(G)$ coincides with its standard linearization if and only if $G$ is Berge-acyclic. This in turn implies that if $G$ is Berge-acyclic, then
$\MP(G)$ is defined by $|V|+(r+2)|E|$ inequalities in the original space.
In~\cite{dPKha18SIOPT}, the authors prove that $\MP(G)$ coincides with its flower relaxation if and only if $G$ is $\gamma$-acyclic. This result implies that if $G$
is $\gamma$-acyclic, then $\MP(G)$ has 
a polynomial-size extended formulation with
at most $|V|+2|E|$ variables and at most $|V|+(r+2)|E|$ inequalities.
%
Subsequently, in~\cite{dPKha23MPA}, the authors 
present a polynomial-size extended formulation for the multilinear polytope of $\beta$-acyclic hypergraphs  with at most $(r-1)|V|+|E|$ variables and at most $(3r-4)|V|+4|E|$ inequalities.

On the other hand, in~\cite{dPDiG23ALG}, the authors prove that \cref{prob MO} is strongly NP-hard over $\alpha$-acyclic hypergraphs.
This result implies that, unless P = NP, one cannot construct, in polynomial time, a polynomial-size extended formulation for the multilinear polytope of $\alpha$-acyclic hypergraphs. 
However, as we detail next, by making further assumptions on the rank of $\alpha$-acyclic hypergraphs, 
one can construct a polynomial-size extended formulation for the multilinear polytope.
In~\cite{WaiJor04,Lau09,BieMun18}, the authors give extended formulations for the convex hull of the feasible set of (possibly constrained) multilinear optimization problems. The size of these extended formulations is parameterized in terms of the ``tree-width'' of their so-called intersection graphs. For the unconstrained case, as detailed in~\cite{dPKha21MOR}, their result can be equivalently stated as follows:
If $G$ is an $\alpha$-acylic hypergraph of rank $r$, with $r = O(\log \poly(|V|, |E|))$, then $\MP(G)$ has a polynomial-size extended formulation, where by $\poly(|V|, |E|)$, we imply a polynomial function in $|V|,|E|$. 
Henceforth, for brevity, whenever for a hypergraph $G=(V, E)$ we have $r = O(\log \poly(|V|, |E|))$, we say that $G$ has \emph{log-poly} rank.

For further results regarding polyhedral relaxations of multilinear sets of degree at least three, 
see~\cite{CraRod16,dPKhaSah20MPC,dPDiG21IJO,dPKha21MOR,Aida22,kim22,dPWal22IPCO,dPWal23MPB}.

\subsection{Pseudo-Boolean optimization}
\label{sec intro Pseudo-Boolean}

We define a \emph{signed hypergraph} $H$ as a pair $(V,S)$, where $V$ is a finite set of nodes and $S$ is a set of \emph{signed edges}.
A \emph{signed edge} $s \in S$ is a pair $(e,\eta_s)$, where $e$ is a subset of $V$ of cardinality at least two, and $\eta_s$ is a map that assigns to each $v \in e$ a \emph{sign} $\eta_s(v) \in \{-1,+1\}$.
The \emph{underlying edge} of a signed edge $s=(e,\eta_s)$ is $e$.
Two signed edges $s=(e,\eta_s)$, $s'=(e',\eta_{s'}) \in S$ are said to be \emph{parallel} if $e=e'$, and they are said to be \emph{identical} if $e=e'$ and $\eta_s=\eta_{s'}$. 
Throughout this paper, we consider signed hypergraphs with no identical signed edges. However, our signed hypergraphs often contain parallel signed edges. 

With any signed hypergraph $H = (V,S)$, and cost vector $c \in \R^{V \cup S}$, we associate the following \emph{pseudo-Boolean optimization problem:}
\begin{align}
\label[problem]{prob PBO}
\tag{$\BPO_\mathrm{pB}$}
\begin{split}
\max & \qquad \sum_{s \in S} c_s \prod_{v \in s} \sigma_s(z_v) \\
{\rm s.t.} & \qquad z \in \{0,1\}^V,
\end{split}
\end{align}
where 
$$
\sigma_s(z_v) := 
\begin{cases}
z_v & \text{ if } \eta_s(v) = +1 \\ 
1-z_v & \text{ if } \eta_s(v) = -1.
\end{cases}
$$
A variety of important applications such as maximum satisfiability problems~\cite{goemans94}, and inference in graphical models~\cite{jordan04} can naturally be formulated as pseudo-Boolean optimization problems.
Problem~\ref{prob PBO} has been extensively studied in the literature; see~\cite{boros02} for a detailed survey of main results. 
Some of the main topics considered in these works are
quadratic pseudo-Boolean optimization problems~\cite{hamhansim84,borCraHam90}, quadratization of general pseudo-Boolean optimization problems~\cite{borcra20}, and special problem classes such as super-modular pseudo-Boolean optimization problems~\cite{crama89}.

As before, we linearize the objective function of~\cref{prob PBO} by introducing a new variable $z_s$, for each signed edge $s \in S$, to obtain an equivalent reformulation of \cref{prob PBO} in a lifted space:
\begin{align}
\label[problem]{prob lPBO}
\tag{$\ell\BPO_\mathrm{pB}$}
\begin{split}
\max & \qquad \sum_{s \in S} {c_s z_s} \\
{\rm s.t.} & \qquad z_s = \prod_{v\in s} {\sigma_s(z_v)} \qquad \forall s \in S \\
& \qquad z \in \{0,1\}^{V \cup S}.
\end{split}
\end{align}
In this paper, we introduce the~\emph{pseudo-Boolean set} 
of the signed hypergraph $H=(V,S)$, as the feasible region of \cref{prob lPBO}:
\begin{equation*} 
\PBS(H) := 
\Big\{ z \in \{0,1\}^{V \cup S} : z_s = \prod_{v\in s} {\sigma_s(z_v)}, \; \forall s \in S \Big\},
\end{equation*}
and we refer to its convex hull as the \emph{pseudo-Boolean polytope} and denote it by $\PBP(H)$. With the objective of constructing strong linear programming (LP) relaxations for pseudo-Boolean optimization problems, in this paper, we study the facial structure of the pseudo-Boolean polytope.

It is important to note that, unlike the multilinear set/polytope, the pseudo-Boolean set/polytope may not be full-dimensional.
For example, if $H=(V,S)$ contains three signed edges 
$s_1=(e,\eta_{s_1})$, $s_2=(e,\eta_{s_2})$, $s_3=(e',\eta_{s_3})$, where $e'= e \setminus \{\bar v\}$ for some $\bar v \in V$, $\eta_{s_1}(v)=\eta_{s_2}(v)=\eta_{s_3}(v)$ for all $v \in e \setminus \{\bar v\}$, and $\eta_{s_1}(\bar v)=-\eta_{s_2}(\bar v)$,
then $\PBS(H)$ and $\PBP(H)$ are not full-dimensional as we have 
$z_{s_1} + z_{s_2} = z_{s_3}$. 
While the pseudo-Boolean polytope is significantly more complex than the multilinear polytope, as we demonstrate in this paper, by studying the facial structure of $\PBP(H)$, one can obtain polynomial-size extended formulations for many instances for which such a formulation  is not known for the corresponding multilinear polytope. In Section~\ref{sec vs} we explore in detail the connections between the two polytopes.


\subsection{Our contributions}

In this paper, we demonstrate how the hypergraph framework pioneered in~\cite{dPKha17MOR} for multilinear optimization is relevant and applicable to pseudo-Boolean optimization as well.
Namely, using signed hypergraphs to represent pseudo-Boolean sets, we present sufficient conditions under which the pseudo-Boolean polytope admits a polynomial-size extended formulation.
Our results unify and extend all prior results on polynomial-size representability of the multilinear polytope~\cite{WaiJor04,Lau09,dPKha18SIOPT,BieMun18,dPKha21MOR,dPKha23MPA}.

We introduce a new technique, which we refer to as the ``recursive inflate-and-decompose'' framework to construct polynomial-size extended formulations for the pseudo-Boolean polytope. Our proposed framework relies on a recursive application of three key ingredients, each of which is of independent interest:
\begin{enumerate}
    \item A sufficient condition for decomposability of pseudo-Boolean polytopes (see~\cref{th decomp}). This is the first result on decomposability of pseudo-Boolean polytopes and serves as a significant generalization of Theorem~4 in~\cite{dPKha23MPA} (Section~\ref{sec: decomp}).
    \item A polynomial-size extended formulation for the pseudo-Boolean polytope of pointed signed hypergraphs (see~\cref{pointedHull}). The pseudo-Boolean polytope of pointed signed hypergraphs is the building block of our extended formulations, which appears as a result of applying our decomposition technique of Part~1 (Section~\ref{sec: pointed}).
    
    \item An operation, which we refer to as ``inflation of signed edges'' (see~\cref{lem inflation}) that we use to transform a large class of signed hypergraphs to those for which our results of Parts~1 and~2 are applicable (Section~\ref{sec: inflate}).
\end{enumerate}

As we detailed in \cref{sec intro multi}, at the time of this writing, the most general sufficient conditions under which one can obtain a polynomial-size extended formulation for the multilinear polytope $\MP(G)$ are:
\begin{itemize}
\item [(i)] $G$ is a $\beta$-acyclic hypergraph~\cite{dPKha23MPA}, \item [(ii)] $G$ is an $\alpha$-acyclic hypergraph with log-poly rank~\cite{Lau09}.
\end{itemize}
It is important to remark that neither of the above sufficient conditions implies the other one. Furthermore, the two results in~\cite{dPKha23MPA} and~\cite{Lau09} have been proven using entirely different techniques. In Section~\ref{sec: framework}, we show that our recursive inflate-and-decompose framework implies as special cases both sufficient conditions~(i) and~(ii) above and extends to many more cases of interest. Below, we summarize these results.

Consider a signed hypergraph $H=(V,S)$.
We define the~\emph{underlying hypergraph} of $H$ as the hypergraph obtained from $H$ by ignoring signs and dropping parallel edges.
In Section~\ref{sec: beta} we prove that, if the underlying hypergraph of $H$ is $\beta$-acyclic, then $\PBP(H)$ has a polynomial-size extended formulation (see Theorem~\ref{extended}). This is a significant generalization of case~(i) above, as it only requires the $\beta$-acyclicity 
of the underlying hypergraph of $H$.
In Section~\ref{sec: alfarank} we prove that if the underlying hypergraph of $H$ is $\alpha$-acyclic and has log-poly rank, then $\PBP(H)$ has a polynomial-size extended formulation (see Theorem~\ref{alphaStuff general}). This result essentially coincides with case~(ii) above. 
In Section~\ref{sec: smallgap}, we introduce the notion of ``gap'' for hypergraphs, which roughly speaking, indicates if it is possible to inflate signed edges in an efficient manner. 
We then show that for signed hypergraphs, if certain gaps are not too large,  by combining results of Theorems~\ref{lem inflation} and~\ref{extended}, one can obtain polynomial-size extended formulations for the pseudo-Boolean polytope (see Propositions~\ref{cor2} and~\ref{cor3}). Finally, in Section~\ref{sec: general} we outline some generalizations and directions of future research. 

It is important to remark that the proofs of all our results regarding the existence of polynomial-size extended formulations for $\PBP(H)$ are constructive and the proposed extended formulations can be constructed in polynomially many operations in $|V|,|S|$.

We would like to conclude this section by further emphasizing on the power of these extended formulations: Not only they serve as polynomial-size LP formulations for special classes of~\cref{prob MO} and~\cref{prob PBO}, they can also be used to construct strong LP \emph{relaxations} for general mixed-integer nonlinear optimization problems whose factorable reformulations contain pseudo-Boolean sets (see for example~\cite{sahTaw03,bel09,misFlo14,IdaNick18,viger18,dPKhaSah20MPC}).

\subsection{Notations and preliminaries}
\label{sec: notation}

In the following, we present all hypergraph terminology and notation that we will use throughout this paper.

\paragraph{Hypergraphs.}

Let $G=(V,E)$ be a hypergraph.
We define the hypergraph obtained from $G$ by \emph{removing} a node $v \in V$ as the hypergraph $G-v$ with set of nodes $V \setminus \{v\}$ and set of edges $\{e - v : e \in E, \ |e-v| \ge 1\}$. 

A node $v \in V$ is a \emph{$\beta$-leaf} of $G$ if the set of the edges of $G$ containing $v$ is totally ordered with respect to inclusion.
A \emph{sequence of $\beta$-leaves} of length $t$ for some  $1 \leq t\leq |V|$ of $G$ is an ordering $v_1, \dots, v_t$ of $t$ distinct nodes of $G$, such that $v_1$ is a $\beta$-leaf of $G$, $v_2$ is a $\beta$-leaf of $G - v_1$, and so on, until $v_t$ is a $\beta$-leaf of $G - v_1 - \dots - v_{t-1}$. 
The hypergraph $G$ is said to be \emph{$\beta$-acyclic} if it has a sequence of $\beta$-leaves of length $|V|$.
An equivalent definition of $\beta$-acyclic hypergraphs can be obtained using the concept of $\beta$-cycle. 
A \emph{$\beta$-cycle} of length $q$ for some $q \ge 3$ in $G$ is a sequence $v_1, e_1, v_2, e_2, \dots , v_q, e_q, v_1$ such that $v_1, v_2, \dots , v_q$ are distinct nodes, $e_1, e_2, \dots , e_q$ are distinct edges, and $v_i$ belongs to $e_{i-1}, e_i$ and no other edge among $e_1, e_2, \dots , e_q$, for all $i = 1, \dots , q$, where $e_0 := e_q$.
A hypergraph is $\beta$-acyclic if and only if it does not contain any $\beta$-cycles \cite{Dur12}. 

A node $v \in V$ is an \emph{$\alpha$-leaf} of $G$ if the set of edges of $G$ containing $v$ has a maximal element for inclusion.
A \emph{sequence of $\alpha$-leaves} of length $t$ for some  $1 \le t\le |V|$ of $G$ is an ordering $v_1, \dots, v_t$ of $t$ distinct nodes of $G$, such that $v_1$ is an $\alpha$-leaf of $G$, $v_2$ is an $\alpha$-leaf of $G - v_1$, and so on, until $v_t$ is an $\alpha$-leaf of $G - v_1 - \dots - v_{t-1}$.
The hypergraph $G$ is said to be \emph{$\alpha$-acyclic} if it has a sequence of $\alpha$-leaves of length $|V|$.
An equivalent definition of $\alpha$-acyclic hypergraphs can be obtained using the concept of $\alpha$-cycles introduced in \cite{JegNdi09}.

We say that $G$ is \emph{connected} if for every $u,w \in V$,
there exists a sequence $v_1, e_1, v_2, e_2, \dots , v_q, e_q, v_{q+1}$ such that $v_1=u$, $v_{q+1}=w$, $v_2, v_3, \dots , v_q \in V$, $e_1, e_2, \dots , e_q \in E$, 
and $e_i$ contains $v_i, v_{i+1}$, 
for all $i = 1, \dots , q$.
%
The \emph{connected components} of $G$ are its maximal connected partial hypergraphs.
A hypergraph $G' = (V',E')$ is a \emph{partial hypergraph} of $G = (V,E)$ if $V' \subseteq V$ and $E' \subseteq E$.

\paragraph{Signed hypergraphs.}
Let $H=(V,S)$ be a signed hypergraph.
In the following, we define some useful operations on signed edges. Let $s=(e,\eta_s) \in S$. 
With a slight abuse of notation, we use set theoretical notation on $s$, with the understanding that $s$ should be replaced with $e$.
For example, we denote by $|s|$ the number $|e|$, we say that $s$ is nonempty if $e$ is nonempty, we write $v \in s$ meaning $v \in e$, and for $U \subseteq V$, we write $s \subseteq U$ (resp., $s \supseteq U$, $s=U$) meaning $e \subseteq U$ (resp., $e \supseteq U$, $e=U$).
Similarly, we might write $s \subseteq s'$, for $s=(e,\eta_s), s'=(e',\eta_{s'}) \in S$, instead of $e \subseteq e'$.
If $v \in s$, we denote by $s - v$ the signed edge $s'=(e',\eta_{s'})$, where $e':=e \setminus \{v\}$, and $\eta_{s'}$ is the restriction of $\eta_s$ that assigns to each $v \in e'$ the sign $\eta_{s'}(v) = \eta_s(v)$.
If $v \notin s$, we denote by $s + v^+$ the signed edge $s'=(e',\eta_{s'})$, where $e':=e \cup \{v\}$, and $\eta_{s'}$ is the extension of $\eta_s$ that assigns to each $u \in e'$ the sign $\eta_{s'}(u) = \eta_s(u)$ and assigns to $v$ the sign $\eta_{s'}(v) = +1$.
Similarly, $s + v^-$ is defined as $s + v^+$ but with $\eta_{s'}(v) = -1$.

We define the signed hypergraph obtained from $H$ by \emph{removing} a node $v \in V$ as the signed hypergraph $H - v$ with set of nodes $V \setminus \{v\}$ and set of signed edges $\{s - v : s \in S, \ |s-v| \ge 1\}$. 

\subsection{Organization}
The remainder of the paper is structured as follows. In Section~\ref{sec vs} we detail on the connections between multilinear optimization and pseudo-Boolean optimization. In Section~\ref{sec: decomp} we present a sufficient condition for decomposability of the pseudo-Boolean polytope. In Section~\ref{sec: pointed} we obtain a polynomial-size extended formulation for the pseudo-Boolean polytope of pointed signed hypergraphs. In Section~\ref{sec: inflate} we introduce the inflation operation. Combining the results of Sections~\ref{sec: decomp}-\ref{sec: inflate}, in Section~\ref{sec: framework}, we introduce our new ``recursive inflate-and-decompose'' framework, which in turn enables us to present polynomial-size extended formulations for the pseudo-Boolean polytope of a large class of signed hypergraphs. Finally, Section~\ref{sec: appendix} is the appendix and contains a technical proof omitted from Section~\ref{sec: pointed}.


\section{Multilinear optimization versus pseudo-Boolean optimization}
\label{sec vs}

In this section we compare the two encodings of binary polynomial optimization introduced in Section~\ref{sec: intro}; \ie multilinear optimization defined by~\cref{prob MO} and pseudo-Boolean optimization defined by \cref{prob PBO}. We also detail on the relationship between 
the multilinear set (resp. multilinear polytope) and the pseudo-Boolean set (resp. pseudo-Boolean polytope).

First, it is simple to see that \cref{prob MO} is a special case of~\cref{prob PBO} obtained by letting $\eta_s(v) = +1$ for every $s \in S$ and $v \in s$.
Similarly, the multilinear set $\MS(G)$ and the multilinear polytope $\MP(G)$ are special cases of the pseudo-Boolean set $\PBS(H)$ and of the pseudo-Boolean polytope $\PBP(H)$, respectively.
Vice versa, it is also possible to go from the pseudo-Boolean setting to the multilinear setting.
In the remainder of the section we discuss this direction, which requires more work and may result in an exponential blowup of the size of the objects considered.


\paragraph{Reformulating \cref{prob PBO} as \cref{prob MO}.}
Each product term in the objective function of \cref{prob PBO} can be expanded and written as:
\begin{align}
\label{eq mapping simple}
\prod_{v\in s} {\sigma_s(z_v)} = \sum_{e \in E} d_e \prod_{v \in e} z_v + d_0,
\end{align}
where the coefficients $d_0$ and $d_e$, for $e \in E$, are $0, \pm 1$ and can be easily computed by expanding the product on the left-hand side of~\eqref{eq mapping simple}.
To rewrite \eqref{eq mapping simple} with explicit coefficients, let $p(s):=\bra{v \in s : \eta_s(v)=+1}$ and, for $i=0,1,2,\dots,|s|$, let $N_i(s):=\bra{t \subseteq s : \eta_s(v)=-1 \ \forall v \in t, \ |t|=i}$.
We then obtain
\begin{align}
\label{eq mapping full}
\prod_{v\in s} {\sigma_s(z_v)} 
= 
\sum_{i=0}^{|s|} \sum_{t \in N_i(s)}
(-1)^{|t|}
\prod_{v \in p(s) \cup t} z_v.
\end{align}
%
It then follows that 
\cref{prob PBO} over a signed hypergraph $H=(V,S)$ can be reformulated (up to a constant in the objective) as \cref{prob MO} over a hypergraph, which we call the ``multilinear hypergraph'' of $H$ and that we define next.

Given a signed hypergraph $H=(V,S)$, the \emph{multilinear hypergraph} of $H$, denoted by $\mh(H)$, is the hypergraph $(V,E)$, where $E$ is constructed as follows: for each $s \in S$, and every $t \subseteq s$ with $\eta_s(v)=-1$ for all $v \in t$, the set $E$ contains $\{v \in s : \eta_s(v)=+1\} \cup t$ if it has cardinality at least two.
It then follows from the definition of $\mh(H)$ that all sets $p(s) \cup t$ that appear in \eqref{eq mapping full} are edges of $\mh(H)$.
We summarize the above discussion in the next observation:

\begin{observation}
\cref{prob PBO} over a signed hypergraph $H=(V,S)$ can be reformulated (up to a constant in the objective) as \cref{prob MO} over $\mh(H)$.
\end{observation}

It is important to note that the number of edges in the definition of $\mh(H)$ that arise from a single signed edge $s \in S$ can be \emph{exponential} in the number of nodes contained in $s$.
%
%
In fact, the reformulation of \cref{prob PBO} as \cref{prob MO} described above
may lead to an exponential increase in the size of the problem. 
One easy example of this is given by the following special case of \cref{prob PBO}:
\begin{align*}
\max & \qquad \prod_{v \in V} (1-z_v) \\
{\rm s.t.} & \qquad z \in \{0,1\}^V.
\end{align*}
In this problem, $S$ contains only one signed edge $s=(V,\eta_s)$ with $\eta_s(v)=-1$ for all $v \in V$.
Once reformulated as \cref{prob MO}, the resulting objective function contains all possible $2^{|V|}-1$ monomials.
However, one can observe that the change of variable $z'_v := 1-z_v$, for all $v \in V$, directly transforms this special case of \cref{prob PBO} to \cref{prob MO} with only one monomial in the objective function, namely $\prod_{v \in V} z'_v$.

Next, we provide a more involved example that is insensitive 
to changes of variables of the form $z'_v := 1-z_v$, for some $v \in V$.
Consider a signed hypergraph $H=(V,S)$, where $S$ contains the two signed edges $s=(V, \eta_s)$ and $t=(V, \eta_t)$ such that $\eta_s(v)=-\eta_t(v)$ for every $v \in V$.
We then consider the special case of \cref{prob PBO} given by
\begin{align*}
\max & \qquad 2 \prod_{v \in s} \sigma_s(z_v) + \prod_{v \in t} \sigma_t(z_v) \\
{\rm s.t.} & \qquad z \in \{0,1\}^V.
\end{align*}
It is then simple to check that, once reformulated as \cref{prob MO}, the number of monomials in the resulting objective function is at least $\sum_{i=1}^m \binom{m}{i} = 2^m-1$, where $m:=\lceil |V|/2 \rceil$.
This bound holds even if we apply any change of variables of the form $z'_v := 1-z_v$, for $v \in U \subseteq V$.

\paragraph{The pseudo-Boolean set as a projection of the multilinear set.}

The pseudo-Boolean set $\PBS(H)$ of a signed hypergraph $H=(V,S)$ can be obtained as a projection of the multilinear set of $\mh(H)$.
Let $\mh(H) = (V,E)$.
In the remainder of the section, we denote by $x_v$, for $v \in V$, the variables corresponding to nodes of $H$, which coincide with nodes of $\mh(H)$. 
These variables are used both for the pseudo-Boolean set and the multilinear set.
We denote the variables of the multilinear set corresponding to edges in $E$ by $y_e$, for $e \in E$.
The variables of the pseudo-Boolean set corresponding to signed edges in $S$ are instead denoted by $z_s$, for $s \in S$.

To obtain the pseudo-Boolean set as a projection of the multilinear set, we need to construct the $z$ variables from the $x,y$ variables.
This can be done by linearizing equality \eqref{eq mapping full}, for every $s \in S.$
Specifically, in \eqref{eq mapping full}, the left-hand side is replaced with $z_s$, and each product $\prod_{v \in p(s) \cup t} z_v$ on the right-hand side is replaced with $y_{p(s) \cup t}$ if $|p(s) \cup t| \ge 2$, and with $x_{p(s) \cup t}$ if $|p(s) \cup t| = 1$.
We denote by $\I$ this system of linear equalities obtained by linearizing equality \eqref{eq mapping full}, for every $s \in S.$
The following projection result is then easy to check:
\begin{observation}
The pseudo-Boolean set $\PBS(H)$ is the projection onto the space of the $x,z$ variables of the sets of points $(x,y,z) \in \R^{V \cup E \cup S}$ satisfying the system of linear equalities $\I$ and such that $(x,y)$ lies in the multilinear set of $\mh(H)$.
\end{observation}



\paragraph{The pseudo-Boolean polytope as a projection of the multilinear polytope.}

The pseudo-Boolean polytope $\PBP(H)$ of a signed hypergraph $H=(V,S)$ can be obtained as a projection of the multilinear polytope of $\mh(H)$.

\begin{observation}
The pseudo-Boolean polytope $\PBP(H)$ is the projection onto the space of the $x,z$ variables of the set of points $(x,y,z) \in \R^{V \cup E \cup S}$ satisfying the system of linear equalities $\I$ and such that $(x,y)$ lies in the multilinear polytope of $\mh(H)$.
\end{observation}



\begin{prf}
Denote by $Q$ the set of points $(x,y,z) \in \R^{V \cup E \cup S}$ satisfying the system of linear equalities $\I$ and such that $(x,y)$ lies in the multilinear polytope of $\mh(H)$.
Denote by $Q'$ the projection of $Q$ onto the space of variables $x,z$.
The inequalities defining the multilinear polytope of $\mh(H)$ only involve variables $x,y$, while the linear equalities in $\I$ involve variables $x,y,z$.
In particular, we obtain that $Q$ is a polytope, and therefore $Q'$ is a polytope too.
Let $(\bar x, \bar z)$ be a vertex of $Q'$.
It suffices to show that $(\bar x, \bar z)$ is in the pseudo-Boolean set $\PBS(H)$.
By definition of $Q'$, there exists $\bar y$ such that $(\bar x, \bar y, \bar z)$ is a vertex of $Q$.
Therefore, there exist $|V|+|E|+|S|$ linearly independent constraints among the linear inequalities defining $Q$ such that $(\bar x, \bar y, \bar z)$ is the unique vector that satisfies these constraints at equality.
Note that the constraints in $\I$ are the only ones containing variables $z$, each such constraint contains exactly one $z$ variable with nonzero coefficient, and there is precisely one constraint for each $z$ variable. 
Thus, $|S|$ of the linearly independent constraints defining $(\bar x, \bar y, \bar z)$ must be the $|S|$ constraints in $\I$.
The remaining $|V|+|E|$ linearly independent constraints defining $(\bar x, \bar y, \bar z)$ must be inequalities defining the multilinear polytope of $\mh(H)$, and these constraints only involve $x,y$ variables.
Hence, $(\bar x, \bar y)$ is a vertex of the multilinear polytope of $\mh(H)$, and hence it is in the multilinear set of $\mh(H)$.
It then follows from our previous discussion that $(\bar x, \bar z)$ is in the pseudo-Boolean set $\PBS(H)$.
\end{prf}





\section{Decomposability of pseudo-Boolean polytopes}
\label{sec: decomp}

In this section, we present a sufficient condition for decomposability of pseudo-Boolean polytopes that we will use to obtain our extended formulations. Our decomposition result is the first known sufficient condition for decomposability of the pseudo-Boolean polytope.
Existing decomposability results for the multilinear polytope are Theorem~1 in~\cite{dPKha18MPA}, 
Theorem~5 in~\cite{dPKha18SIOPT},
Theorem~1 in~\cite{dPKha21MOR},
Theorem~4 in~\cite{dPDiG21IJO}, and
Theorem~4 in~\cite{dPKha23MPA}. Our decomposition result serves as a significant generalization of Theorem~4 in~\cite{dPKha23MPA}.

Consider a signed hypergraph $H = (V,S)$, let $V_1,V_2 \subseteq V$ such that $V = V_1 \cup V_2$, let $S_1 \subseteq \{s \in S : s \subseteq V_1\}$, $S_2 \subseteq \{s \in S : s \subseteq V_2\}$ such that $S = S_1 \cup S_2$.
Let $H_1 := (V_1,S_1)$ and $H_2 := (V_2,S_2)$.
In order to define the concept of decomposability 
we now explain how we can write a vector $z \in \R^{V \cup S}$ by partitioning its components into three sub-vectors.
The vector $z_\cap$ contains the components of $z$ corresponding to nodes and signed edges that are both in $H_1$ and in $H_2$.
The vector $z_1$ contains the components of $z$ corresponding to nodes and signed edges in $H_1$ but not in $H_2$.
Finally, the vector $z_2$ contains the components of $z$ corresponding to nodes and signed edges in $H_2$ but not in $H_1$.
Using these definitions, we can now write a vector $z$ in the space defined by $H$ as $z=(z_1,z_\cap,z_2)$, up to reordering variables.
Similarly, we can write a vector $z$ in the space defined by $H_1$ as $(z_1,z_\cap)$, and a vector $z$ in the space defined by $H_2$ as $z=(z_\cap,z_2)$.
We now define the concept of decomposability.

\begin{figure}[h]
\begin{center}
\includegraphics[width=.8\textwidth]{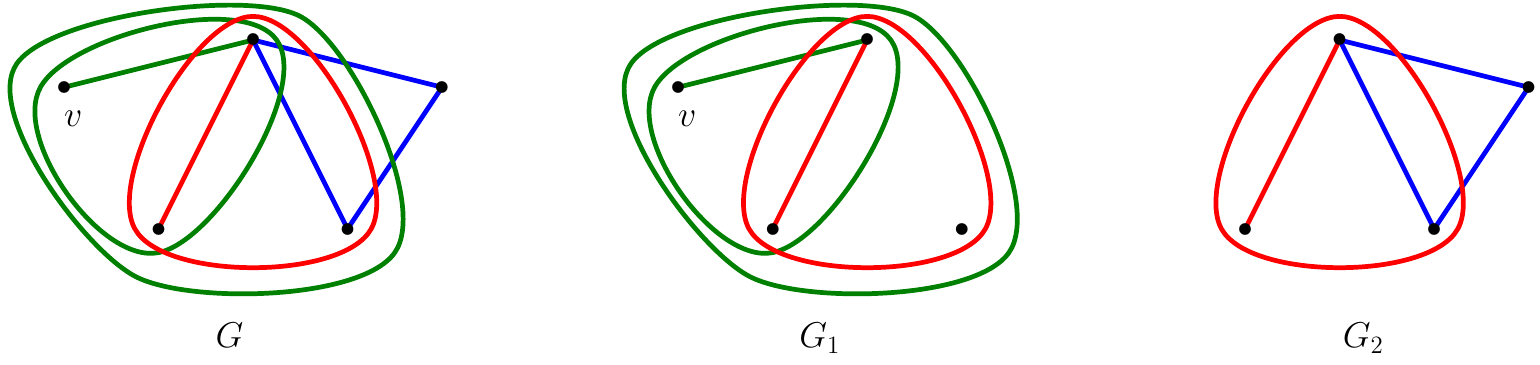}
\caption{Example of signed hypergraphs $H,H_1,H_2$ that follow the construction of \cref{th decomp}. The figure depicts only the underlying hypergraphs of $H,H_1,H_2$, denoted by $G,G_1,G_2$, respectively. The decomposition result holds for any possible signing of the edges of $H$. 
Note that if in the edge of cardinality three containing $v$, all nodes have sign $-1$, then $v$ is not a nest point in the multilinear hypergraph of $H$, and Theorem~4 in~\cite{dPKha23MPA} cannot be applied.}
\label{fig decomp}
\end{center}
\end{figure}

\begin{definition}[Decomposability]
\label{def decomposability}
We say that the pseudo-Boolean polytope $\PBP(H)$ is \emph{decomposable into the pseudo-Boolean polytopes $\PBP(H_1)$ and $\PBP(H_2)$} if the following relation holds
\begin{align} 
\label{eq def decomposability}
(z^1,z^\cap) \in \PBP(H_1) , \ (z^\cap,z^2) \in \PBP(H_2)
\quad \Rightarrow \quad 
(z^1,z^\cap,z^2) \in \PBP(H).
\end{align}
\end{definition}
We observe that the opposite implication in \eqref{eq def decomposability}, from right to left, always holds.
Furthermore, if $\PBP(H)$ is decomposable into $\PBP(H_1)$ and $\PBP(H_2)$, then the system comprised of a description of $\PBP(H_1)$ and a description of $\PBP(H_2)$, is a description of $\PBP(H)$.
%
We are now ready to present our decomposability result, in \cref{th decomp} below.
An example of signed hypergraphs $H,H_1,H_2$ that follow the construction of \cref{th decomp} is given in \cref{fig decomp}.
Throughout this paper, for an integer $k$, we define $[k]:= \{1,\dots,k\}$.

\begin{theorem}
\label{th decomp}
Let $H = (V,S)$ be a signed hypergraph, and assume that the underlying hypergraph of $H$ has a $\beta$-leaf $v$.
Let $s_1 \subseteq s_2 \subseteq \cdots \subseteq s_k$ be the signed edges of $H$ containing $v$, and assume that $S$ contains the signed edges $s_i - v$ such that $|s_i - v| \ge 2$, for every $i \in [k]$.
Then $\PBP(H)$ is decomposable into $\PBP(H_1)$ and $\PBP(H_2)$, where $H_1$ and $H_2$ are defined as follows.
$H_1 := (V_1, S_v \cup P_v)$, where $V_1$ is the underlying edge of $s_k$, $S_v := \{s_1,\dots,s_k\}$, $P_v := \{s_i - v : |s_i - v| \ge 2, \ i \in [k]\}$, and $H_2 := H - v$.
\end{theorem}

\begin{prf}
We assume $k \ge 1$, as otherwise the result is obvious.
Let $z_1$, $z_\cap$, and $z_2$ be defined as above, so that we can write a vector $z$ in the space defined by $H$ as $z=(z_1,z_\cap,z_2)$, a vector $z$ in the space defined by $H_1$ as $(z_1,z_\cap)$, and a vector $z$ in the space defined by $H_2$ as $z=(z_\cap,z_2)$.
In particular, the vector $z_\cap$ contains the components of $z$ corresponding to nodes in $V_1 \setminus \{v\}$ and signed edges in $P_v$, while the vector $z_1$ contains the components of $z$ corresponding to node $v$ and signed edges in $S_v$.

Let $\tilde z$ be such that $(\tilde z_1,\tilde z_\cap) \in \PBP(H_1)$ and $(\tilde z_\cap,\tilde z_2) \in \PBP(H_2)$.
We show $\tilde z \in \PBP(H)$.
To do so, we will write $\tilde z$ explicitly as a convex combinations of vectors in $\PBS(H)$.
In the next claim, we show how a vector in $\PBS(H_1)$ and a vector in $\PBS(H_2)$ can be combined to obtain a vector in $\PBS(H)$.

Define $p_i := s_i - v$ for every $i \in [k]$.
Note that, for some $i=1,\dots,k$, we might have $|p_i|=1$, and so $p_i \notin P_v$, meaning that $p_i$ is not a signed edge of $H$.
In these cases, we have that $p_i$ contains a single node, say $u$, and we call $p_i$ a \emph{signed loop.}
Consistently with our notation for signed edges, we will write 
$z_{p_i}$ and $\sigma_{p_i} (z_{u})$ to denote $z_u$ if $\eta_{p_i}(u)=+1$, or to denote $1-z_u$ if $\eta_{p_i}(u)=-1$.

\begin{claim}
\label{claim combination}
Let $(z_1,z_\cap) \in \PBS(H_1)$ and $(z'_\cap,z'_2) \in \PBS(H_2)$ be such that 
$z_{p_i} = z'_{p_i}$ for every $i \in [k]$.
Then, $(z_1,z'_\cap,z'_2) \in \PBS(H)$.
\end{claim}

\begin{cpf}
We have $(z_1,z'_\cap,z'_2) \in \PBS(H)$ if and only if $(z_1,z'_\cap) \in \PBS(H_1)$ and $(z'_\cap,z'_2) \in \PBS(H_2)$, thus it suffices to show $(z_1,z'_\cap) \in \PBS(H_1)$.
The signed edges of $H_1$ are $S_v \cup P_v$.
Consider first a signed edge $p_i \in P_v$.
In $(z_1,z'_\cap)$, the component corresponding to the signed edge $p_i$ is in $z'_\cap$, and the components corresponding to the nodes in $p_i$ are in $z'_\cap$, thus we need to show the equality $z'_{p_i} = \prod_{u \in p_i} \sigma_{p_i}(z'_u)$, which follows directly from $(z'_\cap,z'_2) \in \PBS(H_2)$.
Consider now a signed edge $s_i \in S_v$.
In $(z_1,z'_\cap)$, the component corresponding to the signed edge $s_i$ is in $z_1$, and the components corresponding to the nodes in $s_i$ are in $z'_\cap$, except for the component corresponding to node $v$ that is in $z_1$.
Thus, we need to show the equality 
$$
z_{s_i} = 
\sigma_{s_i}(z_v) \prod_{u \in s_i \setminus \{v\}} \sigma_{s_i}(z'_u).
$$
Since $s_i$ is a signed edge of $H_1$ and $(z_1,z_\cap) \in \PBS(H_1)$, we know
\begin{align*}
z_{s_i} 
& = \sigma_{s_i}(z_v) \prod_{u \in s_i\setminus\{v\}} \sigma_{s_i}(z_u).
\end{align*}
We then obtain
\begin{align*}
\prod_{u \in s_i\setminus\{v\}} \sigma_{s_i}(z_u)
= \prod_{u \in p_i} \sigma_{p_i}(z_u) 
= z_{p_i} 
= z'_{p_i} 
= \prod_{u \in p_i} \sigma_{p_i}(z'_u) 
= \prod_{u \in s_i\setminus\{v\}} \sigma_{s_i}(z'_u).
\end{align*}
Here, 
the first equality holds by definition of $p_i$ for $i \in [k]$;
the second equality holds since $p_i$ is a signed edge of $H_1$ or a signed loop containing a node of $H_1$;
in the third equality we use the assumption $z_{p_i} = z'_{p_i}$;
the fourth equality holds because $p_i$ is a signed edge of $H_2$ or a signed loop containing a node of $H_2$, and $(z'_\cap,z'_2) \in \PBS(H_2)$;
the last equality holds by definition of $p_i$ for $i \in [k]$.
\end{cpf}

In the remainder of the proof, we show how to write explicitly $\tilde z$ as a convex combination of the vectors in $\PBS(H)$ obtained in \cref{claim combination}.
By assumption, the vector $(\tilde z_1,\tilde z_\cap)$ is in $\conv \PBS(H_1)$. Thus, it can be written as a convex combination of points in $\PBS(H_1)$;
\ie there exists $\mu \ge 0$ such that
\begin{align}
\nonumber
\sum_{(z_1 ,z_\cap) \in \PBS(H_1)} \mu_{(z_1,z_\cap)} & = 1 \\
\label{eq mu conv}
\sum_{(z_1 ,z_\cap) \in \PBS(H_1)} \mu_{(z_1,z_\cap)} (z_1,z_\cap) & = (\tilde z_1,\tilde z_\cap).
\end{align}
Similarly, the vector $(\tilde z_\cap,\tilde z_2)$ is in $\conv \PBS(H_2)$ and it can be written as a convex combination of points in $\PBS(H_2)$;
\ie there exists $\nu \ge 0$ such that 
\begin{align}
\nonumber
\sum_{(z_\cap, z_2) \in \PBS(H_2)} \nu_{(z_\cap,z_2)} & = 1 \\
\label{eq nu conv}
\sum_{(z_\cap, z_2) \in \PBS(H_2)} \nu_{(z_\cap,z_2)} (z_\cap,z_2) & = (\tilde z_\cap,\tilde z_2).
\end{align}

In the remainder of the proof, given $z_{p_1}, \dots, z_{p_k} \in \{0,1\}$, we will consider the number $\max\{j \in [k] : z_{p_j} = 1\} \in \{0,\dots,k\}$, with the understanding that this number equals $0$ when $z_{p_1} = \cdots = z_{p_k} = 0$.
In the next technical claim, we study the sums of the multipliers $\mu$ and $\nu$ corresponding to binary vectors $z$ with a fixed $\max\{j \in [k] : z_{p_j} = 1\}$. 
To do so, we need to introduce the rather technical notation $\ext(i)$, for $i=0,\dots,k$.
%
For $j,t \in \{1,\dots,k\}$ with $j<t$, we say that $p_t$ is an \emph{extension} of $p_j$ if the common nodes of $p_j,p_t$ have the same sign in $p_j$ and in $p_t$; formally, for every $v \in p_j$ we have $\eta_{p_j}(v)= \eta_{p_t}(v)$.
We can now define $\ext(0)$ as the set of indices $t \in \{1,\dots,k\}$, such that there is no $j \in \{1,\dots,t-1\}$ such that $p_t$ is an extension of $p_j$.
For $i=1,\dots,k$, we define $\ext(i)$ as the set of indices $t \in \{i+1,\dots,k\}$, such that $p_t$ is a \emph{minimal extension} of $p_i$, that is, $p_t$ is an extension of $p_i$, and there is no $j \in \{i+1,\dots,t-1\}$ such that $p_t$ is an extension of $p_j$.
Note that, if instead there exists $j \in \{i+1,\dots,t-1\}$ such that $p_t$ is an extension of $p_j$, then $p_j$ is also an extension of $p_i$.
Observe that $\ext(k) = \emptyset$.
%
%
To familiarize with the definition, note that in the special case where all signed edges $s_1,\dots,s_k$ have all nodes with a $+1$ sign, i.e., $\eta_{s_i}(v)=+1$ for every $v \in s_i$ and $i=1,\dots,k$, we have 
$\ext(0)=\bra{1}$ and $\ext(i)=\bra{i+1}$, for $i=1,\dots,k$.
We are now ready to present our technical claim.
\begin{claim}
\label{claim sum mu nu}
For $i \in \{0,\dots,k\}$, we have
\begin{align*}
\sum_{\substack{(z_1 ,z_\cap) \in \PBS(H_1) \\ \max\{j \in [k] : z_{p_j} = 1\} = i}} \mu_{(z_1,z_\cap)} 
=
\sum_{\substack{(z_\cap, z_2) \in \PBS(H_2) \\ \max\{j \in [k] : z_{p_j} = 1\} = i}} \nu_{(z_\cap, z_2)} 
=
\begin{cases}
1 - \sum_{t \in \ext(0)} \tilde z_{p_{t}} & \qquad \text{if $i = 0$} \\
\tilde z_{p_i} - \sum_{t \in \ext(i)} \tilde z_{p_{t}} & \qquad \text{if $i \in \{1,\dots,k\}$}.
\end{cases}
\end{align*}
\end{claim}

\begin{cpf}
By considering the component of \eqref{eq mu conv} corresponding to $p_i$, for $i \in [k]$, we obtain
\begin{align*}
\sum_{\substack{(z_1 ,z_\cap) \in \PBS(H_1) \\ z_{p_i}=1}} \mu_{(z_1,z_\cap)} 
= \tilde z_{p_i}.
\end{align*}

We first consider the case $i=k$.
We have
\begin{align*}
\sum_{\substack{(z_1 ,z_\cap) \in \PBS(H_1) \\ \max\{j \in [k] : z_{p_j} = 1\} = k}} \mu_{(z_1,z_\cap)} 
=
\sum_{\substack{(z_1 ,z_\cap) \in \PBS(H_1) \\ z_{p_k}=1}} \mu_{(z_1,z_\cap)} 
= \tilde z_{p_k}
=
\tilde z_{p_k} - \sum_{t \in \ext(k)} \tilde z_{p_t}.
\end{align*}

Next, we consider the case $i \in \{1,\dots,k-1\}$.
We have
\begin{align*}
\sum_{\substack{(z_1 ,z_\cap) \in \PBS(H_1) \\ \max\{j \in [k] : z_{p_j} = 1\} = i}} \mu_{(z_1,z_\cap)} 
& =
\sum_{\substack{(z_1 ,z_\cap) \in \PBS(H_1) \\ z_{p_i}=1 \\ z_{p_{t}}=0 \ \forall t \in \{i+1,\dots,k\}}} \mu_{(z_1,z_\cap)} \\
& =
\sum_{\substack{(z_1 ,z_\cap) \in \PBS(H_1) \\ z_{p_i}=1 \\ z_{p_{t}}=0 \ \forall t \in \ext(i)}} \mu_{(z_1,z_\cap)} \\
& =
\sum_{\substack{(z_1 ,z_\cap) \in \PBS(H_1) \\ z_{p_i}=1}} \mu_{(z_1,z_\cap)} 
-
\sum_{t \in \ext(i)}
\sum_{\substack{(z_1 ,z_\cap) \in \PBS(H_1) \\ z_{p_t}=1}} \mu_{(z_1,z_\cap)} \\
&
= \tilde z_{p_i} - \sum_{t \in \ext(i)} \tilde z_{p_{t}}.
\end{align*}
To check the second equality, we assume $(z_1 ,z_\cap) \in \PBS(H_1)$ with $z_{p_i}=1$ and $z_{p_{t}}=0$ for all $t \in \ext(i)$ and we show that $z_{p_{t}}=0$ for all $t \in \bra{i+1,\dots,k}$.
For a contradiction, assume that there exists $t \in \bra{i+1,\dots,k}$ with $z_{p_{t}}=1$.
Assume that $t$ is the smallest index of this type.
Then $t \notin \ext(i)$.
Now $p_t$ must be an extension of $p_i$, as otherwise $z_{p_i}=1$ implies $z_{p_t}=0$, a contradiction.
Hence, $p_t$ is not a minimal extension of $p_i$, thus there exists $j \in \{i+1,\dots,t-1\}$ such that $p_t$ is an extension of $p_j$.
By the choice of $t$ we have $z_{p_{j}}=0$, which implies $z_{p_{t}}=0$, since $p_t$ is an extension of $p_j$. 
We have obtained a contradiction.

To check the third equality, we assume $(z_1 ,z_\cap) \in \PBS(H_1)$ with $z_{p_i}=1$ and $z_{p_{t}}=1$ for some $t \in \ext(i)$ and we show that $z_{p_{t}}=1$ for precisely one index $t \in \ext(i)$.
For a contradiction, assume there exist two distinct $t,j \in \ext(i)$ with $z_{p_{t}}=z_{p_{j}}=1$.
Without loss of generality, assume $j<t$.
By definition of $\ext(i)$, both $p_j$ and $p_t$ are extensions of $p_i$, and $p_t$ is not an extension of $p_j$.
But then $z_{p_j}=1$ implies $z_{p_t}=0$, a contradiction.

Lastly, we consider the case $i=0$.
We have
\begin{align*}
\sum_{\substack{(z_1 ,z_\cap) \in \PBS(H_1) \\ \max\{j \in [k] : z_{p_j} = 1\} = 0}} \mu_{(z_1,z_\cap)} 
&=
\sum_{\substack{(z_1 ,z_\cap) \in \PBS(H_1) \\ z_{p_t}=0 \ \forall t \in \{1,\dots,k\}}} \mu_{(z_1,z_\cap)} \\
&=
\sum_{\substack{(z_1 ,z_\cap) \in \PBS(H_1) \\ z_{p_t}=0 \ \forall t \in \ext(0)}} \mu_{(z_1,z_\cap)} \\
&=
\sum_{\substack{(z_1 ,z_\cap) \in \PBS(H_1)}} \mu_{(z_1,z_\cap)} - \sum_{t \in \ext(0)} \sum_{\substack{(z_1 ,z_\cap) \in \PBS(H_1) \\ z_{p_t}=1}} \mu_{(z_1,z_\cap)} \\
&=
1 - \sum_{t \in \ext(0)} \tilde z_{p_{t}}.
\end{align*}
The equalities above follow with arguments similar to those in the previous case.

The statement for $\nu$ follows symmetrically, starting with \eqref{eq nu conv} rather than \eqref{eq mu conv}.
\end{cpf}

For ease of notation, we define, for $i \in \{0,\dots,k\}$, the quantity $m(i)$ to be the sum obtained in \cref{claim sum mu nu}, i.e.,
\begin{align*}
m(i) := 
\begin{cases}
1 - \sum_{t \in \ext(0)} \tilde z_{p_{t}} & \qquad \text{if $i = 0$} \\
\tilde z_{p_i} - \sum_{t \in \ext(i)} \tilde z_{p_{t}} & \qquad \text{if $i \in \{1,\dots,k\}$}.
\end{cases}
\end{align*}

We are now ready to define the multipliers $\lambda$ that will allow us to write $\tilde z$ as a convex combination of vectors in $\PBS(H)$ obtained in \cref{claim combination}.
For every $(z_1,z_\cap) \in \PBS(H_1)$ and $(z'_\cap,z'_2) \in \PBS(H_2)$ such that
$z_{p_j} = z'_{p_j}$ for every $j \in [k]$, 
we define 
\begin{align*}
\lambda_{(z_1,z'_\cap,z'_2)} :=
\begin{cases}
\frac{\mu_{(z_1,z_\cap)} \nu_{(z'_\cap,z'_2)}}
{m(i)} & \qquad \text{if $m(i) \neq 0$} \\
 0 & \qquad \text{if $m(i) = 0$},
\end{cases}
\end{align*}
where $i := \max\{j \in [k] : z_{p_j} = 1\} = \max\{j \in [k] : z'_{p_j} = 1\}$.
%
%
In the next claim, we show that the multipliers $\lambda$ are nonnegative and sum to one.

\begin{claim}
\label{claim sum one}
We have $\lambda \ge 0$ and 
\begin{align*}
\sum_{\substack{(z_1 ,z_\cap) \in \PBS(H_1) \\ (z'_\cap, z'_2) \in \PBS(H_2) \\ z_{p_j} = z'_{p_j} \ \forall j \in [k]}} 
\lambda_{(z_1,z'_\cap,z'_2)} =1.
\end{align*}
\end{claim}

\begin{cpf}
It follows from \cref{claim sum mu nu} that $m(i) \ge 0$ for all $i \in \{0,\dots,k\}$.
Thus, using the definition of $\lambda$, we obtain $\lambda \ge 0$.
%
Using \cref{claim sum mu nu}, we obtain
\begin{align*}
\sum_{\substack{(z_1 ,z_\cap) \in \PBS(H_1) \\ (z'_\cap, z'_2) \in \PBS(H_2) \\ z_{p_j} = z'_{p_j} \ \forall j \in [k]}} 
\lambda_{(z_1,z'_\cap,z'_2)} 
& = 
\sum_{i=0,\dots,k} \ 
\sum_{\substack{(z_1 ,z_\cap) \in \PBS(H_1) \\ \max\{j \in [k] : z_{p_j} = 1\} = i}} \
\sum_{\substack{(z'_\cap, z'_2) \in \PBS(H_2) \\ \max\{j \in [k] : z'_{p_j} = 1\} = i}} \
\lambda_{(z_1,z'_\cap,z'_2)} \\
& = 
\sum_{\substack{i=0,\dots,k \\ m(i) \neq 0}} \ 
\sum_{\substack{(z_1 ,z_\cap) \in \PBS(H_1) \\ \max\{j \in [k] : z_{p_j} = 1\} = i}} \
\sum_{\substack{(z'_\cap, z'_2) \in \PBS(H_2) \\ \max\{j \in [k] : z'_{p_j} = 1\} = i}} \
\lambda_{(z_1,z'_\cap,z'_2)} \\
& = 
\sum_{\substack{i=0,\dots,k \\ m(i) \neq 0}} \ 
\sum_{\substack{(z_1 ,z_\cap) \in \PBS(H_1) \\ \max\{j \in [k] : z_{p_j} = 1\} = i}} \
\sum_{\substack{(z'_\cap, z'_2) \in \PBS(H_2) \\ \max\{j \in [k] : z'_{p_j} = 1\} = i}} \
\frac{\mu_{(z_1,z_\cap)} \nu_{(z'_\cap,z'_2)}}{m(i)} \\
& = 
\sum_{\substack{i=0,\dots,k \\ m(i) \neq 0}} \ 
\frac{1}{m(i)}
\pare{
\sum_{\substack{(z_1 ,z_\cap) \in \PBS(H_1) \\ \max\{j \in [k] : z_{p_j} = 1\} = i}} 
\mu_{(z_1,z_\cap)}}
\pare{
\sum_{\substack{(z'_\cap, z'_2) \in \PBS(H_2) \\ \max\{j \in [k] : z'_{p_j} = 1\} = i}} 
\nu_{(z'_\cap,z'_2)} } \\
& = 
\sum_{\substack{i=0,\dots,k \\ m(i) \neq 0}} \ 
\frac{(m(i))^2}{m(i)} 
= 
\sum_{\substack{i=0,\dots,k \\ m(i) \neq 0}} \ 
m(i) 
= 
\sum_{i=0,\dots,k} \ 
m(i) 
= 1.
\end{align*}
The last equality $\sum_{i=0,\dots,k} m(i) = 1$ can be seen using the definition of $m(i)$, because each index in $\{1,\dots,k\}$ is exactly in one set among $\ext(0),\dots,\ext(k)$.
\end{cpf}

Our last claim, which concludes the proof of the theorem, shows that the multipliers $\lambda$ yield $\tilde z$ as a convex combination of the vectors in $\PBS(H)$ obtained in \cref{claim combination}.

\begin{claim}
\label{claim conv comb}
We have 
\begin{align}
\label{eq to verify}
(\tilde z_1,\tilde z_\cap,\tilde z_2) = 
\sum_{\substack{(z_1 ,z_\cap) \in \PBS(H_1) \\ (z'_\cap, z'_2) \in \PBS(H_2) \\ z_{p_j} = z'_{p_j} \ \forall j \in [k]}} 
\lambda_{(z_1,z'_\cap,z'_2)} 
(z_1,z'_\cap,z'_2),
\end{align}
\end{claim}

\begin{cpf}
Using the definition of $\lambda$, we rewrite \eqref{eq to verify} in the form
\begin{align*}
(\tilde z_1,\tilde z_\cap,\tilde z_2)
& = 
\sum_{\substack{(z_1 ,z_\cap) \in \PBS(H_1) \\ (z'_\cap, z'_2) \in \PBS(H_2) \\ z_{p_j} = z'_{p_j} \ \forall j \in [k]}} 
\lambda_{(z_1,z'_\cap,z'_2)} 
(z_1,z'_\cap,z'_2) \\
& =
\sum_{i=0,\dots,k} \ 
\sum_{\substack{(z_1 ,z_\cap) \in \PBS(H_1) \\ \max\{j \in [k] : z_{p_j} = 1\} = i}} \
\sum_{\substack{(z'_\cap, z'_2) \in \PBS(H_2) \\ \max\{j \in [k] : z'_{p_j} = 1\} = i}} \
\lambda_{(z_1,z'_\cap,z'_2)} (z_1,z'_\cap,z'_2) \\
& =
\sum_{\substack{i=0,\dots,k \\ m(i) \neq 0}} \ 
\sum_{\substack{(z_1 ,z_\cap) \in \PBS(H_1) \\ \max\{j \in [k] : z_{p_j} = 1\} = i}} \
\sum_{\substack{(z'_\cap, z'_2) \in \PBS(H_2) \\ \max\{j \in [k] : z'_{p_j} = 1\} = i}} \
\lambda_{(z_1,z'_\cap,z'_2)} (z_1,z'_\cap,z'_2) \\
& = 
\sum_{\substack{i=0,\dots,k \\ m(i) \neq 0}} \ 
\sum_{\substack{(z_1 ,z_\cap) \in \PBS(H_1) \\ \max\{j \in [k] : z_{p_j} = 1\} = i}} \
\sum_{\substack{(z'_\cap, z'_2) \in \PBS(H_2) \\ \max\{j \in [k] : z'_{p_j} = 1\} = i}} \
\frac{\mu_{(z_1,z_\cap)} \nu_{(z'_\cap,z'_2)}}{m(i)} (z_1,z'_\cap,z'_2).
\end{align*}
We now verify the obtained equality, first for components $\tilde z_1$, and then for components $\tilde z_\cap,\tilde z_2$.
We start with components $\tilde z_1$.
We obtain
\begin{align*}
\tilde z_1 
& = 
\sum_{\substack{i=0,\dots,k \\ m(i) \neq 0}} \ 
\sum_{\substack{(z_1 ,z_\cap) \in \PBS(H_1) \\ \max\{j \in [k] : z_{p_j} = 1\} = i}} \
\sum_{\substack{(z'_\cap, z'_2) \in \PBS(H_2) \\ \max\{j \in [k] : z'_{p_j} = 1\} = i}} \
\frac{\mu_{(z_1,z_\cap)} \nu_{(z'_\cap,z'_2)}}{m(i)} z_1\\
& = 
\sum_{\substack{i=0,\dots,k \\ m(i) \neq 0}} \ 
\frac{1}{m(i)}
\pare{
\sum_{\substack{(z_1 ,z_\cap) \in \PBS(H_1) \\ \max\{j \in [k] : z_{p_j} = 1\} = i}} 
\mu_{(z_1,z_\cap)}z_1}
\pare{
\sum_{\substack{(z'_\cap, z'_2) \in \PBS(H_2) \\ \max\{j \in [k] : z'_{p_j} = 1\} = i}} 
\nu_{(z'_\cap,z'_2)}} \\
& = 
\sum_{\substack{i=0,\dots,k \\ m(i) \neq 0}} \ 
\sum_{\substack{(z_1 ,z_\cap) \in \PBS(H_1) \\ \max\{j \in [k] : z_{p_j} = 1\} = i}} 
\mu_{(z_1,z_\cap)}z_1 \\
& = 
\sum_{i=0,\dots,k} \ 
\sum_{\substack{(z_1 ,z_\cap) \in \PBS(H_1) \\ \max\{j \in [k] : z_{p_j} = 1\} = i}} 
\mu_{(z_1,z_\cap)}z_1 \\
& = 
\sum_{(z_1 ,z_\cap) \in \PBS(H_1)} 
\mu_{(z_1,z_\cap)} z_1.
\end{align*}
Here, the third equality holds due to \cref{claim sum mu nu}.
The forth equality holds because, from \cref{claim sum mu nu}, if $m(i)=0$, then all $\mu_{(z_1,z_\cap)}$ in the second sum equal zero.
The obtained equality is implied by \eqref{eq mu conv}.

Next, we consider components $\tilde z_\cap,\tilde z_2$.
We obtain
\begin{align*}
(\tilde z_\cap,\tilde z_2) 
& = 
\sum_{\substack{i=0,\dots,k \\ m(i) \neq 0}} \ 
\sum_{\substack{(z_1 ,z_\cap) \in \PBS(H_1) \\ \max\{j \in [k] : z_{p_j} = 1\} = i}} \
\sum_{\substack{(z'_\cap, z'_2) \in \PBS(H_2) \\ \max\{j \in [k] : z'_{p_j} = 1\} = i}} \
\frac{\mu_{(z_1,z_\cap)} \nu_{(z'_\cap,z'_2)}}{m(i)} (z'_\cap,z'_2)\\
& = 
\sum_{\substack{i=0,\dots,k \\ m(i) \neq 0}} \ 
\frac{1}{m(i)}
\pare{
\sum_{\substack{(z_1 ,z_\cap) \in \PBS(H_1) \\ \max\{j \in [k] : z_{p_j} = 1\} = i}} 
\mu_{(z_1,z_\cap)}}
\pare{
\sum_{\substack{(z'_\cap, z'_2) \in \PBS(H_2) \\ \max\{j \in [k] : z'_{p_j} = 1\} = i}} 
\nu_{(z'_\cap,z'_2)} (z'_\cap,z'_2)} \\
& = 
\sum_{\substack{i=0,\dots,k \\ m(i) \neq 0}} \ 
\sum_{\substack{(z'_\cap, z'_2) \in \PBS(H_2) \\ \max\{j \in [k] : z'_{p_j} = 1\} = i}} 
\nu_{(z'_\cap,z'_2)} (z'_\cap,z'_2) \\
& = 
\sum_{i=0,\dots,k} \ 
\sum_{\substack{(z'_\cap, z'_2) \in \PBS(H_2) \\ \max\{j \in [k] : z'_{p_j} = 1\} = i}} 
\nu_{(z'_\cap,z'_2)} (z'_\cap,z'_2) \\
& = 
\sum_{(z'_\cap, z'_2) \in \PBS(H_2)} 
\nu_{(z'_\cap,z'_2)} (z'_\cap,z'_2).
\end{align*}
Here, the third equality holds due to \cref{claim sum mu nu}.
The forth equality holds because, from \cref{claim sum mu nu}, if $m(i)=0$, then all $\nu_{(z'_\cap,z'_2)}$ in the second sum equal zero.
The obtained equality is implied by \eqref{eq nu conv}.
\end{cpf}

%
We should remark that the overall structure of the proof of \cref{th decomp} is similar to that of Theorem~4 in~\cite{dPKha23MPA}.
The difference lies in the key construction of the proof, that is, the multipliers introduced to construct the convex combination.
Due to the presence of the signed edges, these new multipliers are significantly more involved than the ones presented in the proof of Theorem~4 in~\cite{dPKha23MPA}.
In the special case where all signs are positive, the new multipliers simplify to the ones in the proof of Theorem~4 in~\cite{dPKha23MPA}.

\section{Pointed signed hypergraphs}
\label{sec: pointed}

The signed hypergraph $H_1$ defined in the statement of Theorem~\ref{th decomp} plays a key role in our convex hull characterizations. In this section, we prove that the pseudo-Boolean polytope $\PBP(H_1)$ has a polynomial-size extended formulation. Together with the decomposition technique of Theorem~\ref{th decomp}, this result enables us to obtain polynomial-size extended formulations for the pseudo-Boolean polytope of a large family of signed hypergraphs.

\begin{definition}[Pointed signed hypergraph]
\label{def pointed}
Consider a signed hypergraph $H=(V,S)$ and let $v \in V$ be a $\beta$-leaf of the underlying hypergraph of $H$. Denote by $S_v$ the set of all signed edges in $S$ containing $v$. 
Define $P_v:=\{s-v: s \in S_v, |s|\geq 3\}$.
We say  that $H$ is \emph{pointed} at $v$ (or is a \emph{pointed signed hypergraph}) if $V$ coincides with the underlying edge of the signed edge of maximum cardinality in $S_v$ and $S = S_v \cup P_v$.
%
\end{definition}

From this definition it follows that the 
signed hypergraph $H_1$ in the statement of Theorem~\ref{th decomp} is pointed at $v$ (the hypergraph $G_1$ in Figure~\ref{fig decomp} is an example of the underlying hypergraph of a pointed signed hypergraph).
In order to characterize the pseudo-Boolean polytope of a signed hypergraph $H$ pointed at $v$, we first characterize the pseudo-Boolean polytope of a simpler type of signed hypergraphs corresponding to faces of $\PBP(H)$ defined by $z_v = 0$ or $z_v =1$. An extended formulation for $\PBP(H)$ can then be obtained using disjunctive programming.


\subsection{The pseudo-Boolean polytope of nested signed hypergraphs}

In this section, we characterize the pseudo-Boolean polytope of nested signed hypergraphs in the original space. This result will then enable us to obtain a polynomial-size extended formulation for the pseudo-Boolean polytope of pointed signed hypergraphs.

\begin{definition}[Nested signed hypergraph]
\label{def nested}
Let $H=(V,S)$ be a signed hypergraph with $V=\{v_1, \dots, v_n\}$.
Denote by $E$ the set of underlying edges of $S$ (note that several signed edges may collapse to the same underlying edge). 
Define $e_k := \{v_1, \dots v_{k+1}\}$ for all $k \in [n-1]$ and $\bar E:= \{e_1, \dots, e_{n-1}\}$.
We say that $H$ is a \emph{nested signed hypergraph} 
if it satisfies the following conditions:
\begin{itemize}
\item [(N1)] $E = \bar E$,
\item  [(N2)] for any signed edge $s \in S$ with underlying edge $e_k$ for some $k \in [n-1]$, the following two signed edges are also present in $S$:
$\ell(s)$, obtained from $s$ by flipping the sign of $v_{k+1}$,
and
$p(s):= s - v_{k+1}$. 
\end{itemize}
\end{definition}
%
\begin{remark}\label{rem1}
Nested signed hypergraphs have two important properties that we will use to obtain our extended formulations: 
\begin{enumerate}
    \item As we will prove in Proposition~\ref{nestedHull}, the pseudo-Boolean polytope of a nested signed hypergraph $H=(V,S)$ in the original space is defined by
    $\frac{|S|}{2}+1$ equalities, and $|S|+2|V|-4$ inequalities; \ie $2(|S|+|V|-1)$ inequalities.
    \item Let $H=(V,S)$ be a signed hypergraph such that for any $s,s' \in S$ we either have $s \subseteq s'$ or $s' \subseteq s$. It then follows that by adding at most $2|S|(|V|-2)+|S|$ signed edges to $H$ we obtain a nested signed hypergraph. To see this, note that to satisfy property~(N2) of nested signed hypergraphs, for each edge $s \in S$ we need to add at most $2(|s|-2)+1$ signed edges to $H$. Together with case~1 above, this in turn implies that $\PBP(H)$ has a polynomial-size extended formulation with at most $2|S|(|V|-1)+|V|$ variables and at most $4 |S|(|V|-1)+2|V|$ inequalities.
\end{enumerate}
\end{remark}

We now proceed with characterizing the pseudo-Boolean polytope of nested signed hypergraphs. To this end, we first introduce some notation.
We denote by $\E_k$, $k \in [n-1]$, the set consisting of signed edges in $S$ whose underlying edges are $e_k$. 
For each $k \in [n-1] \setminus \{1\}$, we define 
$$\E^{+}_k :=\Big\{s \in \E_k: \eta_s(v_{k+1})=+1\Big\},
\qquad \text{and} \qquad
\E^{-}_k:=\Big\{s \in \E_k: \eta_s(v_{k+1})=-1\Big\}.$$ 
For any $s \in S$, we define 
$$N(s):= \{s' \in S: s =p(s')\},$$
where as before for any $s' \in S$ with underlying edge $e_k$, we define $p(s')=s'-v_{k+1}$.
Notice that by property~(N2) of nested signed hypergraphs, $|N(s)|$ equals zero or two.  

The following proposition characterizes the pseudo-Boolean polytope of nested signed hypergraphs in the original space. 
Due to its length and technical nature, we have deferred the proof to the appendix; \ie Section~\ref{sec: appendix}. While our proof relies on a standard disjunctive programming technique~\cite{Bal98} followed by a projection step using Fourier-Motzkin elimination, the novelty of the proof lies in the manner Fourier-Motzkin elimination is implemented. Namely, it is well-understood that a generic application of Fourier-Motzkin elimination leads to a rapid increase in the number of inequalities defining the polyhedron. In our proof, 
the auxiliary variables are projected out in a specific order so that the projection does not contain redundant inequalities.

%
\begin{proposition}\label{nestedHull}
Let $H=(V,S)$ with $n:=|V|$ be a nested signed hypergraph. Suppose that $\E_1 = \{q_1, q_2, q_3, q_4\}$  with $\eta_{q_1}(v_1)=\eta_{q_1}(v_2)= \eta_{q_2}(v_1)=\eta_{q_3}(v_2)= +1$ and
$\eta_{q_2}(v_2)=\eta_{q_3}(v_1)= \eta_{q_4}(v_1)=\eta_{q_4}(v_2)= -1$. Then
the pseudo-Boolean polytope $\PBP(H)$ is given by: 
    \begin{align}
    &z_{s} + z_{\ell(s)} = z_{p(s)}, \quad \forall s \in \E^{+}_k, k \in [n-1]\setminus\{1\}\label{eqk1}\\
    &z_s \geq 0, \quad \forall s \in S\label{eqk4}\\  
    &     \sum_{s \in \E^{+}_k}{z_s}\leq z_{v_{k+1}}, \quad \forall k \in [n-1] \setminus \{1\}\label{n1}\\
     &    \sum_{s \in \E^{-}_k}{z_s}\leq 1-z_{v_{k+1}},  \quad \forall k \in [n-1] \setminus \{1\}\label{n2}\\
     &      z_{q_1} + z_{q_2} = z_{v_1}, \; 
    z_{q_1} + z_{q_3} = z_{v_2}, \; 
    z_{q_3} + z_{q_4} = 1-z_{v_1}.\label{fl}
    \end{align}
\end{proposition}

\begin{remark}\label{rem2}
In Proposition~\ref{nestedHull}, the assumption on the structure of $\E_1$ is not restrictive. Recall that by property~(N2) of a nested signed hypergraph we must either have $|\E_1| =2$ or $|\E_1|=4$. Moreover, if $|\E_1|= 2$, then we must either have $\E_1=\{q_1,q_3\}$ or $\E_1=\{q_2, q_4\}$. Without loss of generality, suppose that 
$\E_1=\{q_1,q_3\}$. Then to obtain a description of $\PBP(H)$ in the original space, it suffices to project out $z_{q_2}, z_{q_4}$ from system~\eqref{eqk1}--\eqref{fl}.  To do so, notice that $z_{q_2}, z_{q_4}$ appear only in the following constraints of system~\eqref{eqk1}--\eqref{fl}:
\begin{align*}
    & z_{q_2} \geq 0, \; z_{q_4} \geq 0 \\
    & z_{q_1}+z_{q_2}=z_{v_1}, \; z_{q_3}+z_{q_4}=1-z_{v_1}.
\end{align*}
Projecting out $z_{q_2}, z_{q_4}$, we obtain:
$$
z_{q_1} \leq z_{v_1}, \; z_{q_3} \leq 1-z_{v_1},
$$
which together with the remaining equalities and inequalities of system~\eqref{eqk1}--\eqref{fl} gives a description of $\PBP(H)$ in the original space.
\end{remark}

It is important to note that, in spite of its simple structure, the constraint matrix of the pseudo-Boolean polytope of a nested signed hypergraph is \emph{not} totally unimodular. The following example demonstrates this fact. 
In this example, for ease of notation, we write a signed edge in a compact form by listing its nodes with their signs as superscripts.
For example, the signed edge $s=(e,\eta_s)$, where $e=\{v_1,v_2,v_3\}$ and $\eta_s(v_1)=+1$, $\eta_s(v_2)=-1$, $\eta_s(v_3)=+1$, will be written compactly as $s=\{v_1^+, v_2^-, v_3^+\}$.

\begin{example}\label{eg1}
Consider the nested signed hypergraph $H=(V,S)$, with $V=\{v_1, v_2,v_3,v_4\}$, $\E_1 =\{s_1,s_2,s_3,s_4\}$, $s_1=\{v_1^{+},v_2^{+}\}$, $s_2=\{v^{-}_1,v^{+}_2\}$, $s_3=\{v_1^+,v_2^-\}$, $s_4=\{v_1^-,v_2^-\}$, $\E_2=\{s_5, s_6, s_7, s_8\}$, 
$s_5=\{v_1^+,v_2^+,v_3^+\}$, $s_6=\{v_1^+,v_2^+,v_3^-\}$
$s_7=\{v_1^-,v_2^+,v_3^-\}$, $s_8=\{v_1^-,v_2^+,v_3^+\}$,
$\E_3=\{s_9,s_{10},s_{11},s_{12}\}$, 
$s_9=\{v_1^+,v_2^+,v_3^-,v_4^+\}$,
$s_{10}=\{v_1^+,v_2^+,v_3^-,v_4^-\}$,
$s_{11}=\{v_1^-,v_2^+,v_3^+,v_4^+\}$,
$s_{12}=\{v_1^-,v_2^+,v_3^+,v_4^-\}$.
Then by Proposition~\ref{nestedHull}, the following constraints are present in the description of $\PBP(H)$:
\begin{align*}
&z_{s_5} + z_{s_6} = z_{s_1}    \\
&z_{s_9} + z_{s_{10}} = z_{s_6} \\    
&z_{s_{11}} + z_{s_{12}} = z_{s_8}\\    
&z_{s_5}+z_{s_8} \leq z_{v_3} \\
&z_{s_9}+z_{s_{11}} \leq  z_{v_4}  
\end{align*}
Consider the submatrix of the above equalities and inequalities corresponding to $z_{s_5},z_{s_6},z_{s_8},z_{s_9},z_{s_{11}}$. It can be checked that the determinant of this submatrix equals 2, implying the constraint matrix of $\PBP(H)$ is not totally unimodular. In addition, each row and each column of this submatrix has two non-zero entries, implying it is not a balanced matrix either~\cite{Cor01b}. 
\end{example}


\subsection{The pseudo-Boolean polytope of pointed signed hypergraphs}

Thanks to the compact description for the pseudo-Boolean polytope of nested signed hypergraphs given by Proposition~\ref{nestedHull}, we now provide a polynomial-size extended formulation for the pseudo-Boolean polytope of pointed signed hypergraphs. Recall that in a signed hypergraph $H=(V,S)$ pointed at $\bar v$ for some $\bar v \in V$, we have $S=S_{\bar v} \cup P_{\bar v}$; the set $S_{\bar v}$ consists of all signed edges containing $\bar v$, and the set $P_{\bar v}$ has the form $P_{\bar v}=\{s-\bar v: s\in S_{\bar v}, |s|\geq 3\}$. 

\begin{remark}\label{rem3}
Consider a pointed signed hypergraph $H=(V,S)$. There are two special cases for which a polynomial-size description of $\PBP(H)$ follows from previously known results:
\begin{itemize}
    \item [(i)] Suppose that for each $v \in V$ we have $\eta_s(v) = \eta_{s'}(v)$
for all signed edges $s,s'$ in $H$ containing $v$. Then (possibly after a one-to-one linear transformation) the explicit description of $\PBP(H)$ in the original space  consisting of at most $5|V|+2$ inequalities is given in Theorem~5 of~\cite{dPKha23MPA}. 

 \item [(ii)] Suppose that all signed edges in $S_{\bar v}$ are parallel. It then follows that for any  $s,s' \in S$  we either have $s \subseteq s'$ or $s' \subseteq s$. Then by Part~2 of Remark~\ref{rem1} and definition of $P_{\bar v}, S_{\bar v}$, the pseudo-Boolean polytope $\PBP(H)$ has a polynomial-size extended formulation with at most $|S|(|V|-1)+|V|$ variables and at most $2(|S|(|V|-1)+|V|)$ inequalities.
\end{itemize}
\end{remark}

The following theorem gives a polynomial-size extended formulation for the pseudo-Boolean polytope of general pointed signed hypergraphs.

\begin{theorem}\label{pointedHull}
Let $H=(V,S)$ be a pointed signed hypergraph. 
Then the pseudo-Boolean polytope $\PBP(H)$ has a polynomial-size extended formulation with at most $2|V|(|S|+1)$ variables and at most $4(|S|(|V|-2)+|V|)$ inequalities. Moreover, all coefficients and right-hand side constants in the system defining 
$\PBP(H)$ are $0,\pm 1$.
\end{theorem}

\begin{prf}
Suppose that $H =(V,S)$ is pointed at $\bar v$ for some $\bar v \in V$, implying that $\bar v$ is a $\beta$-leaf of the underlying hypergraph of $H$ and $S = S_{\bar v} \cup P_{\bar v}$. Define the signed hypergraph 
$L=(V\setminus\{\bar v\}, P_{\bar v})$. By definition of $P_{\bar v}$, for any two edges $s,s' \in P_{\bar v}$ we either have $s \subseteq s'$ or $s'\subseteq s$. Hence after the addition of at most $2|P_{\bar v}|(|V|-3)+|P_{\bar v}|$ edges to $L$
we can construct a nested signed hypergraph denoted by $L'$. Denote by $P'_{\bar v}$ the set of signed edges of $L'$. Note that $|P'_{\bar v}| = 2|P_{\bar v}|(|V|-2) \leq |S|(|V|-2)$. Define $H':=(V, S_{\bar v} \cup P'_{\bar v})$. Since a description of $\PBP(H')$ serves as an extended formulation for $\PBP(H)$, to complete the proof, it suffices to show that $\PBP(H')$ has a polynomial-size extended formulation. Denote by $H^0$ (resp. $H^1$) the signed hypergraph corresponding to the face of $\PBP(H')$ with $z_{\bar v}=0$ (resp. $z_{\bar v}=1$). We then have:
$$
\PBP(H')=\conv\Big(\PBP(H^0) \cup \PBP(H^1)\Big).
$$
Denote by $\bar z$ the vector consisting of $z_v$, $v \in V \setminus \{\bar v\}$, and $z_s$, $s \in P'_{\bar v}$. It then follows that:
\begin{align*}
 \PBP(H^0) = \Big\{& z \in \R^{V \cup S_{\bar v} \cup P'_{\bar v}} : z_{\bar v}=0, \; z_s =0, \forall s \in S_{\bar v}\; {\rm with}\; \eta_s(\bar v)=+1, z_s=z_{p(s)}, \forall s \in S_{\bar v}\; {\rm with} \; \\
& \eta_s(\bar v)=-1,\;\bar z \in \PBP(L'_{\bar v})\Big\}\\
\PBP(H^1) = \Big\{& z \in \R^{V \cup S_{\bar v} \cup P'_{\bar v}} : z_{\bar v}=1, \; z_s =z_{p(s)}, \forall s \in S_{\bar v}\; {\rm with} \;\eta_s(\bar v)=+1, z_s=0, \forall s \in S_{\bar v} \; {\rm with} \;\\
& \eta_s(\bar v)=-1,\; \bar z \in \PBP(L'_{\bar v})\Big\}.
\end{align*}
By Propositon~\ref{nestedHull}, the polytope $\PBP(L'_{\bar v})$ is given by system~\eqref{eqk1}-\eqref{fl} and this description has at most 
$2(|S|(|V|-2)+|V|)$ inequalities. For notational simplicity, let us write
$\PBP(L'_{\bar v})$ compactly as 
$$\PBP(L'_{\bar v})=\{\bar z: A \bar z \leq b, C \bar z = d\}.$$
Then, using Balas' formulation for the union of polytopes~\cite{Bal98}, we obtain a polynomial-size extended formulation for $\PBP(H')$:
\begin{align*}
    \PBP(H') =\Big\{& z \in \R^{V \cup S_{\bar v}\cup P'_{\bar v}}: \; \exists (z^0,z^1,z,\lambda) \; {\rm s.t.}\; z=z^0+z^1, \; \lambda_0+\lambda_1=1, \; 
    z^0_{\bar v}=0, \\
    & z^0_s=0, \; \forall s \in S_{\bar v}\;{\rm with}\; \eta_s(\bar v)=+1,\;
    z^0_s = z^0_{p(s)}, \; \forall s \in S_{\bar v} \;{\rm with}\; \eta_s(\bar v)=-1,\\
    & z^1_{\bar v}=\lambda_1, z^1_s=z^1_{p(s)}, \forall s \in S_{\bar v} \;{\rm with}\; \eta_s(\bar v)=+1,\;
    z^1_s = 0, \; \forall s \in S_{\bar v} \;{\rm with}\; \eta_s(\bar v)=-1,\\
    & Az^0\leq b\lambda_0, \; C z^0=d \lambda_0,  Az^1\leq b\lambda_1, \; C z^1=d \lambda_1, \; \lambda_0, \lambda_1 \geq 0\Big\}.
\end{align*}
The size of the above extended formulation can be further reduced by projecting out variables $\lambda_0,\lambda_1, z^0$ using the equalities $z_{\bar v}= \lambda_1$, $\lambda_0+\lambda_1=1$, and $z=z^0+z^1$. Hence, we obtain an extended formulation for $\PBP(H')$ with at most 
\begin{align}\label{novar}
2(|V|+|S_{\bar v}|+|P'_{\bar v}|)& =  2(|V|+|S_{\bar v}|+2|P_{\bar v}|(|V|-2)) = 2(|V|+|S|+|P_{\bar v}|(2|V|-5))\nonumber\\
&\leq  2(|V|+|S|(|V|-\frac{3}{2})) \leq 2|V|(|S|+1),
\end{align} 
variables and at most
\begin{equation}\label{noeq}
2(2(|S|(|V|-2)+|V|)) \leq 4(|S|(|V|-2)+|V|),  
\end{equation}
inequalities. The second equality in~\eqref{novar} follows from $|S|=|S_{\bar v}|+|P_{\bar v}|$ and the inequality in~\eqref{novar} follows from $|P_{\bar v}| \leq \frac{|S|}{2}$. 

By Proposition~\ref{nestedHull}, all coefficients and right-hand side constants in the system defining $\PBP(L'_{\bar v})$ are $0,\pm 1$. 
This together with the fact that
$\lambda_0 = 1-z_{\bar v}$ and $\lambda_1 = z_{\bar v}$ implies that the same statement holds for the extended formulation of $\PBP(H')$ and as a result for the extended formulation of $\PBP(H)$ as well.
\end{prf}

\section{Inflation of signed edges}
\label{sec: inflate}

Our decomposition result stated in Theorem~\ref{th decomp} relies on the key assumption that the underlying hypergraph of $H$ has at least one $\beta$-leaf. 
We argue that in many cases of interest, this restrictive assumption can be removed.
To this end, in the following, we introduce an operation on signed hypergraphs, which we refer to as the~\emph{inflation of signed edges}.
Starting from a signed hypergraph $H$ whose underlying hypergraph does not contain any $\beta$-leaves, by inflating certain edges, we obtain a new signed hypergraph $H'$ whose underlying hypergraph has a sequence of $\beta$-leaves. By relating the extended formulations of $\PBP(H)$ and $\PBP(H')$, we are able to obtain polynomial-size extended formulations for the pseudo-Boolean polytope of various classes of signed hypergraphs whose underlying hypergraphs \emph{contain} $\beta$-cycles.




\begin{definition}[Inflation operation]
\label{def inflation}
Let $H=(V,S)$ be a signed hypergraph, let $s \in S$, and let $e \subseteq V$ such that $s \subset e$.
Denote by $I(s,e)$ the set of all possible signed edges $s'$ parallel to $e$ such that $\eta_s(v) = \eta_{s'}(v)$ for every $v \in s$.
We say that $H'=(V,S')$ is obtained from $H$ by \emph{inflating} $s$ to $e$ if $S' = S \cup I(s,e) \setminus \{s\}$. We also say that $H'$ is obtained from $H$ via an \emph{inflation operation}. 
\end{definition}

\begin{figure}[h]
\begin{center}
\includegraphics[width=.5\textwidth]{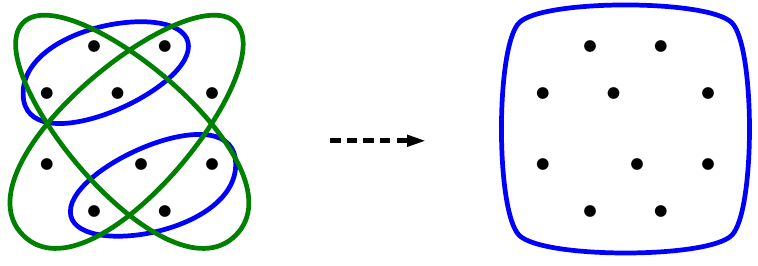}
\caption{The underlying hypergraph of a signed hypergraph $H=(V,S)$ is shown on the left.  The underlying hypergraph of the signed hypergraph $H'$ obtained from $H$ by inflating all signed edges in $S$ to $V$ is shown on the right.}
\label{fig inflation}
\end{center}
\end{figure}

We refer the reader to \cref{fig inflation} for an example of the inflation operation.
The next theorem indicates that if an extended formulation for $\PBP(H')$ is available, one can obtain an extended formulation for $\PBP(H)$ as well.

\begin{theorem}
\label{lem inflation}
Let $H=(V,S)$ be a signed hypergraph, let $s \in S$, and let $e \subseteq V$ such that $s \subset e$.
Let $H'=(V,S')$ be obtained from $H$ by inflating $s$ to $e$.
Then an extended formulation of $\PBP(H)$ can be obtained by juxtaposing an extended formulation of $\PBP(H')$ and the equality constraint
\begin{align}
\label{eq inflation sum}
    z_s = \sum_{s' \in I(s,e)}{z_{s'}}.
\end{align}
Moreover, if $\PBP(H')$ has a polynomial-size extended formulation and $|e|-|s| = O(\log\poly(|V|,|S|))$, then $\PBP(H)$ has a polynomial-size extended formulation as well.
\end{theorem}

\begin{prf}
Given a signed hypergraph $H$ and a set $C$, we denote by $\proj_{H}(C)$ the orthogonal projection of the set $C$ onto the subspace of variables corresponding to nodes and signed edges of $H$. 
Let $A(z,y) \le b$ be an extended formulation of $\PBP(H')$. 
Let $H''=(V,S'')$, where $S'' := S \cup I(s,e) = S' \cup \{s\}$.
It suffices to prove that $A(z,y) \le b$ together with \eqref{eq inflation sum} is an extended formulation of $\PBP(H'')$.
That is, we need to show that 
$$
Q := \proj_{H''} \{(z,y) : A(z,y) \le b, \ \eqref{eq inflation sum}\}
= \PBP(H'').
$$
Since $S \subseteq S''$, it follows that an extended formulation of $\PBP(H'')$ is also an extended formulation of $\PBP(H)$.
For every binary $z$, we have $z \in \PBP(H'')$ if and only if 
\begin{align}
\label{eq inflation sum expanded}
z_s = 
\prod_{v \in s} \eta_s(z_v)
=
\prod_{v \in s} \eta_s(z_v) \prod_{v \in e \setminus s}(z_v + (1-z_v))
= \sum_{s' \in I(s,e)}{z_{s'}}.
\end{align}
Thus, it suffices to show that the vertices of $Q$ are binary.
We have
\begin{align}
\label{eq key to lemma}
\begin{split}
Q
& = \proj_{H''} \{(z,y) : A(z,y) \le b\} \cap \{z : \eqref{eq inflation sum}\} \\
& = \proj_{H'} \{(z,y) : A(z,y) \le b\} \cap \{z : \eqref{eq inflation sum}\} \\
& = \PBP(H') \cap \{z : \eqref{eq inflation sum}\},
\end{split}
\end{align}
where the first equality holds since all variables that appear with nonzero coefficients in \eqref{eq inflation sum} correspond to signed edges in $H''$,
the second equality holds because $\{s\} = S'' \setminus S'$ and so the variable $z_s$ does not appear in the system $A(z,y) \le b$, and
the third equality holds since $A(z,y) \le b$ is an the extended formulation of $\PBP(H')$.


Let $\bar z$ be a vertex of $Q$. 
The variable $\bar z_s$ is only present in $\eqref{eq inflation sum}$, thus we have $\bar z_s = \sum_{s' \in I(s,e)}{\bar z_{s'}}$.
The other constraints of the system defining $Q$ that are active at $\bar z$ are all from the projection of the system $A(z,y) \le b$, thus they determine a vertex of $\PBP(H')$, which is binary.
This implies that all components of $\bar z$, except for $\bar z_s$, are binary.
It then follows from \eqref{eq inflation sum expanded} that $\bar z_s$ is binary as well.
%
%

Finally, since $|I(s,e)| \leq 2^{|e|-|s|}$, we conclude that, $\PBP(H)$ has a polynomial-size extended formulation if $\PBP(H')$ has a polynomial-size extended formulation and $|e|-|s| = O(\log \poly(|V|,|S|))$.  
\end{prf}


The inflation operation is of independent interest as it enables us to obtain polynomial-size extended formulations for the pseudo-Boolean polytope of certain signed hypergraphs. 
The following proposition is an illustration of this fact:

\begin{proposition}\label{cor1}
Consider a signed hypergraph $H=(V,S)$.
Suppose that each $s \in S$ contains at least $|V|-k$ nodes. Then $\PBP(H)$ has an extended formulation with $O(2^{k}|V||S|)$ variables and inequalities.
In particular, if $k= O(\log \poly(|V|, |S|))$, then $\PBP(H)$ has a polynomial-size extended formulation. Moreover, all coefficients and right-hand side constants in the system defining 
$\PBP(H)$ are $0,\pm 1$.
\end{proposition}

\begin{prf} 
Consider the signed hypergraph $H'=(V,S')$ obtained from $H$ by inflating every $s \in S$ with $|s| < |V|$ to $V$. We then have $|S'|\leq 2^k|S|$. By Part~2 of Remark~\ref{rem1},
the polytope $\PBP(H')$ has a polynomial-size extended formulation with at most $2^{k+1}|S|(|V|-1)+|V|$ variables and at most $2^{k+2}|S|(|V|-1)+2|V|$ inequalities. Therefore, from Lemma~\ref{lem inflation} we deduce that an extended formulation for $\PBP(H)$ is obtained by adding at most $|S|$ equality constraints (containing $|S|$ additional variables) to the extended formulation of $\PBP(H')$. Since all coefficients and right-hand side constants in the system defining $\PBP(H')$ as well as in equalities~\eqref{eq inflation sum} are $0,\pm 1$, the statement holds for the system defining $\PBP(H)$.
\end{prf}

Figure~\ref{fig inflation} illustrates the proof idea in Proposition~\ref{cor1}. That is, starting from a signed hypergraph $H$ whose underlying edges have large cardinalities, one can inflate all edges to their union, and subsequently use Part~$(ii)$ of Remark~\ref{rem3} to obtain a polynomial-size extended formulation for the pseudo-Boolean polytope $\PBP(H)$.


\section{The recursive inflate-and-decompose framework}
\label{sec: framework}

With the objective of constructing extended formulations for the pseudo-Boolean polytope,
by combining our decomposition scheme of Theorem~\ref{th decomp}, our convex hull characterization of Theorem~\ref{pointedHull}, and the inflation operation of Theorem~\ref{lem inflation}, we introduce a new framework, which we refer to as the \emph{recursive inflate-and-decompose} framework. We show that this framework enables us to obtain polynomial-size extended formulations for the pseudo-Boolean polytope of large families of signed hypergraphs.

As we detail in this section, our new framework unifies and extends all existing results on polynomial-size extended formulations for the convex hull of the feasible region of Problem~\eqref{prob lMO}~\cite{WaiJor04,Lau09,BieMun18, dPKha18SIOPT,dPKha21MOR,dPKha23MPA}. 
In the following, we present our recursive inflate-and-decompose framework.

\bigskip

\noindent
\textbf{The Recursive inflate-and-decompose (RID) framework}
\vspace{0.05in}

\noindent
\textbf{Input.} 
A signed hypergraph $H=(V,S)$.

\noindent
\textbf{Output.} 
An extended formulation for $\PBP(H)$.

\noindent
\textbf{Step 0.} 
Set $H^{(0)} := H$, $i:=0$.

\noindent
\textbf{Step 1.} 
If we can obtain ${\bar H^{(i)}}$ from ${H^{(i)}}$ via a number of inflation operations, such that a suitable extended formulation for $\PBP({\bar H^{(i)}})$ is available, then by \cref{lem inflation} we are done. 
Otherwise, proceed to Step 2.

\noindent
\textbf{Step 2.} 
Choose a node $\bar v$ of $H^{(i)}$.
If $\bar v$ is a $\beta$-leaf of the underlying hypergraph of $H^{(i)}$, then set ${\bar H^{(i)}} := H^{(i)}$ and proceed to Step~3.
Otherwise, construct  $\bar H^{(i)}$ from ${H^{(i)}}$ via inflation operations, such that $v$ is a $\beta$-leaf of the underlying hypergraph of ${\bar H^{(i)}}$. 
By \cref{lem inflation}, it suffices to find an  extended formulation for $\PBP({\bar H^{(i)}})$.

\noindent
\textbf{Step 3.} 
Using  \cref{th decomp}, decompose $\PBP(\bar H^{(i)})$ into $\PBP(\bar H^{(i)}_1)$ and $\PBP(\bar H^{(i)}_2)$, where $\bar H^{(i)}_1$ denotes the signed hypergraph containing node $\bar v$. 
Since an extended formulation for $\PBP(\bar H^{(i)}_1)$ is given by \cref{pointedHull}, it suffices to find an extended formulation for $\PBP(\bar H^{(i)}_2)$.
Set $H^{(i+1)} := \bar H^{(i)}_2$, increment $i$ by one, and go to Step~1. 

\bigskip

Clearly, Steps~1 and~2 of the RID framework can be performed in many different ways. That is, we have not specified how node $\bar v$ in Step 2 should be selected or how the inflation operations of Step~1 and Step~2 should be performed. A simple way to obtain a $\beta$-leaf in Step 2 is to inflate each signed edge containing $\bar v$ to the union of all signed edges containing $\bar v$.
While the RID framework is fairly general, we are naturally interested in specifying conditions under which this framework results in a polynomial-size extended formulation for the pseudo-Boolean polytope. To this end, first note that at each iteration, one node of the signed hypergraph $H^{(i)}$ is removed. Hence, 
the number of iterations of the RID framework is upper bounded by the number of nodes of $H$. It then follows that RID provides a polynomial-size extended formulation for $\PBP(H)$ if the following conditions are satisfied:
\begin{itemize}
    \item [(I)] In Step~1, the algorithm should terminate, only if a polynomial-size extended formulation for $\PBP(\bar H^{(i)})$
    is available.
    \item [(II)] The total number of new edges introduced as a result of inflation operations in Steps~1 and~2 is upper bounded by a polynomial in $|V|,|S|$.
\end{itemize}

In the remainder of this section, we consider various types of signed hypergraphs for which one can customize the RID framework so that conditions~(I)--(II) above are satisfied and hence polynomial-size extended formulations for the corresponding pseudo-Boolean polytope can be constructed. 




Our results are stated in terms of easily verifiable conditions on the structure of the underlying hypergraphs, namely, $\beta$-acyclic hypergraphs, $\alpha$-acyclic hypergraphs with log-poly ranks, and 
certain classes of hypergraphs with log-poly ``gaps''. In each case, we first provide a self contained proof without any reference to the RID framework. Subsequently, we explain how each proof can be obtained via a specific implementation of RID.

\subsection{$\beta$-acyclic hypergraphs}
\label{sec: beta}
Consider a signed hypergraph $H = (V, S)$. In this section, we show that if the underlying hypergraph of 
$H$ has a sequence of $\beta$-leaves, then $\PBP(H)$ can be described in terms of pseudo-Boolean polytopes of simpler signed hypergraphs. In particular, if the underlying hypergraph of  $H$ is $\beta$-acyclic, then $\PBP(H)$ has an extended formulation of size polynomial in $|V|,|S|$. 
These theorems serve as significant generalizations of the main results in~\cite{dPKha23MPA}.


\begin{theorem}\label{tbetaleaves}
Let $H=(V,S)$ be a signed hypergraph of rank $r$, and let $v_1, \dots, v_t$ for some $t\geq 1$ be a sequence of $\beta$-leaves of the underlying hypergraph of $H$. Then an extended formulation for $\PBP(H)$ is given by a description of $\PBP(H-v_1-\cdots-v_t)$
together with a system of at most $4rt(2|S|+1)$ inequalities consisting of at most 
$2rt(2|S|+1)$ variables. Moreover, all coefficients and right-hand side constants in the latter are $0,\pm 1$.
\end{theorem}

\begin{prf}
Denote by $S_{v_1}$ the set of all signed edges of $H$ containing $v_1$. Since $v_1$ is a $\beta$-leaf of the underlying hypergraph of $H$, the set $S_{v_1}$ is totally ordered.
Define the signed hypergraph $H'_1:=(V,S'_1)$
with $S'_1 = S \cup P_{v_1}$, where $P_{v_1}:=\{s-v_1: s\in S_{v_1}, |s| \geq 3\}$.  
Clearly, an extended formulation for $\PBP(H'_1)$ serves as an extended formulation for  $\PBP(H)$
as well.
Now define the pointed signed hypergraph $H_{v_1}:= (V_1, S_{v_1}\cup P_{v_1})$, where $V_1$ denotes the underlying edge of a signed edge of maximum cardinality in $S_{v_1}$. We then have $H'_1 = H_{v_1} \cup (H-v_1)$, where we used the fact that $H'_1-v_1= H-v_1$. Hence by Theorem~\ref{th decomp}, the pseudo-Boolean polytope $\PBP(H'_1)$ is decomposable into pseudo-Boolean polytopes $\PBP(H_{v_1})$ and $\PBP(H-v_1)$.

Next consider the signed hypergraph $H-v_1$. By definition $v_2$ is a $\beta$-leaf of the underlying hypergraph of $H-v_1$. Denote by $S_{v_2}$ the set of signed edges $H-v_1$ containing $v_2$.
Define $S'_2 := (S \setminus S_{v_1})\cup P_{v_2}$, where $P_{v_2}:=\{s-v_2: s\in S_{v_2}, |s|\geq 3\}$, and define the signed hypergraph $H'_2:=(V\setminus\{v_1\},S'_2)$. Again, an extended formulation for $\PBP(H'_2)$ serves as an extended formulation for
$\PBP(H-v_1)$ as well.
Define the pointed signed hypergraph  $H_{v_2}:= (V_2, S_{v_2}\cup P_{v_2})$, where $V_2$ denotes the underlying edge of a signed edge of maximum cardinality in $S_{v_2}$.
Then by Theorem~\ref{th decomp}, $\PBP(H'_2)$ is decomposable into $\PBP(H_{v_2})$ and $\PBP(H-v_1-v_2)$.

By a recursive application of the above argument after $t$ times, we conclude that an extended formulation for $\PBP(H)$ is given by 
a description of $\PBP(H-v_1-\cdots-v_t)$ together with a system of inequalities defining $\PBP(H_{v_i})$ for all $i \in \{1,\dots, t\}$, where $H_{v_i}$ is a signed hypergraph pointed at $v_i$
defined as $H_{v_i}:=(V_i, S_{v_i}\cup P_{v_i})$, where $S_{v_i}$ denotes the set of signed edge of $H-v_1-\cdots -v_{i-1}$ (we define $H-v_1-v_0 = H$) containing $v_i$, $P_{v_i} :=\{s-v_i: s \in S_{v_i}, |s|\geq 3\}$, and $V_i$ denotes the underlying edge of a signed edge of maximal cardinality in $S_{v_i}$.

By Theorem~\ref{pointedHull} the polytopes $\PBP(H_{v_i})$, $i \in \{1,\dots,t\}$ have polynomial-size
extended formulations each of which consists of at most $2|V_i| (|S_{v_i}|+|P_{v_i}|+1) \leq 2 r(2|S|+1)$ variables and at most $4 (|S_{v_i}|+|P_{v_i}|)(|V_i|-2)+4 |V_i| \leq 8(r-2)|S|+4r \leq 4r(2|S|+1)$ inequalities, where the inequalities follow from $|V_i| \leq r$ and $|P_{v_i}| \leq |S_{v_i}| \leq |S|$.
\end{prf}

\begin{remark}
The proof of Theorem~\ref{tbetaleaves} follows from the RID framework. To see this, first note that the node $\bar v$ chosen in Step~2 of the $i$th iteration of the RID framework is the $\beta$-leaf $v_i$ for all $i \in [t]$. Moreover, at iteration $t+1$, the algorithm terminates at Step~1 with a description of $\PBP(H-v_1-\cdots-v_t)$. In this case, no inflation operation at any step of the RID is needed.
\end{remark}

\begin{theorem}\label{extended}
Let $H=(V,S)$ be a signed hypergraph of rank $r$ whose underlying hypergraph is $\beta$-acyclic. Then the pseudo-Boolean polytope $\PBP(H)$ has a polynomial-size extended formulation with at most $2r(2|S|+1)|V|$ variables, and at most $4r(2|S|+1)|V|$ inequalities.
Moreover, all coefficients and right-hand side constants in the system defining 
$\PBP(H)$ are $0,\pm 1$.
\end{theorem}

\begin{prf}
Since the underlying hypergraph of $H$ is $\beta$-acyclic, it has a sequence of
$\beta$-leaves of length $|V|$.  The proof then follows directly from~\cref{tbetaleaves}. 
\end{prf}

\begin{remark}
Let $G=(V,E)$ be a $\beta$-acyclic hypergraph of rank $r$.
Theorem~1 in~\cite{dPKha23MPA} gives an extended formulation for the multilinear polytope of $G$ with at most $(r-1)|V|$ variables and at most $(3r-4)|V|+4|E|$ inequalities. 
Theorem~\ref{extended} is a significant generalization of this result as it only requires the $\beta$-acyclicity of the underlying hypergraph of $H$. Indeed the multilinear hypergraph $\mh(H)$ may contain many $\beta$-cycles. Nonetheless, this generalization has a cost: while the size of the extended formulation for $\MP(G)$ is quadratic in $|V|, |E|$, the size of the extended formulation for $\PBP(H)$ is cubic in $|V|,|S|$.  
\end{remark}

The next example illustrates the power of~\cref{extended}.

\begin{example}\label{expl3}
    Consider the signed hypergraph $H=(V,S)$ with $V=\{v_1, \ldots, v_n\}$ for some $n\geq 4$, and $S=\{s_1, s_2, s_3, s_4\}$,
    where $s_i=(e_i, \eta_{s_i})$ for $i \in \{1,\ldots,4\}$, $e_1=e_2 = \{v_1, \ldots, v_{n-1}\}$, $e_3=e_4 = \{v_2, \ldots, v_n\}$, 
    $\eta_{s_1}(v) =1$ for $v \in \{v_1, \ldots v_{n-2}\}$,  $\eta_{s_1}(v_{n-1}) =-1$,
    $\eta_{s_2}(v) =1$ for $v \in \{v_1, v_3,\ldots v_{n-1}\}$, $\eta_{s_2}(v_{2}) =-1$, $\eta_{s_3}(v) =1$ for $v \in \{v_2, \ldots v_{n-2}, v_n\}$, $\eta_{s_3}(v_{n-1}) =-1$, $\eta_{s_4}(v) =1$ for $v \in \{v_3, \ldots v_{n}\}$, and $\eta_{s_1}(v_{2}) =-1$. It then follows that the underlying hypergraph of $H$ is $\beta$-acyclic and hence a polynomial-size extended formulation for $\PBP(H)$ is given by Theorem~\ref{extended}.
    It can be checked that the multilinear hypergraph of $H$ is given by $\mh(H)=(V,E)$
    with $E=\{f_1, \ldots, f_6\}$, where $f_1=f_2=V$, $f_3= \{v_1, \ldots, v_{n-2}\}$,
    $f_4=\{v_1, v_3, \ldots, v_{n-1}\}$, $f_5=\{v_2,\ldots, v_{n-2}, v_n\}$, and $f_6=\{v_3,\ldots, v_6\}$. It then follows that 
    $\mh(H)$ contains a $\beta$-cycle of length four given by $C=v_1 f_3 v_2 f_5 v_n f_6 v_{n-1} f_4 v_1$.
    Therefore, Theorem~1 in~\cite{dPKha23MPA} cannot be used to obtain an extended formulation for the multilinear polytope of $\mh(H)$. See Figure~\ref{fig th5} for an illustration.
\end{example}

\begin{figure}[h]
\begin{center}
\includegraphics[width=.8\textwidth]{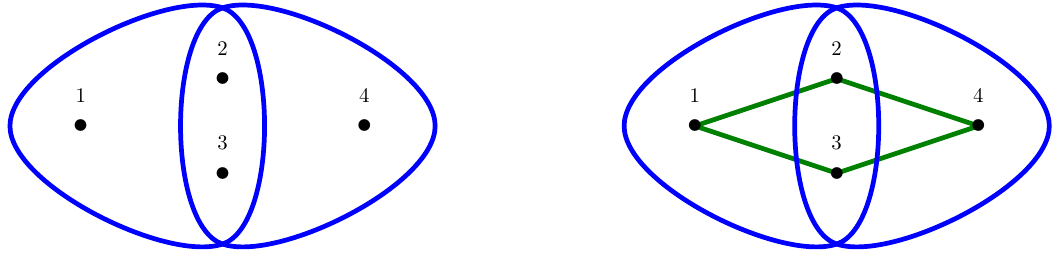}
\caption{The hypergraphs of Example~\ref{expl3} with $n=4$. The underlying hypergraph of the signed hyerpgraph $H$,  which is $\beta$-acyclic is shown on the left. The multilinear hypergraph of $H$ that contains a $\beta$-cycle of length four is shown on the right.}
\label{fig th5}
\end{center}
\end{figure}

\subsection{$\alpha$-acyclic hypergraphs with log-poly ranks}
\label{sec: alfarank}
As we mentioned before, $\alpha$-acyclic hypergraphs are the most general type of acyclic hypergraphs. In~\cite{dPDiG23ALG}, the authors prove 
that unless P = NP, one cannot construct, in polynomial time, a polynomial-size extended formulation for the multilinear polytope of $\alpha$-acyclic hypergraphs.
In~\cite{WaiJor04,Lau09,BieMun18}, the authors give extended formulations for the convex hull of the feasible set of (possibly constrained) multilinear optimization problems. For the unconstrained case, their result can be equivalently stated as follows (see Theorem~5 in~\cite{dPKha21MOR}):
%
If $G=(V, E)$ is an $\alpha$-acylic hypergraph of rank $r$ such that $r = O(\log \poly (|V|,|E|))$, then $\MP(G)$ has a polynomial-size extended formulation. 
In the following we show that this result follows from Theorems~\ref{th decomp}--\ref{lem inflation}; that is, it is a special case of the RID framework.






\begin{theorem}\label{alphaStuff general}
Let $H = (V, S)$ be a signed hypergraph of rank $r$ whose underlying hypergraph is $\alpha$-acyclic.
Then $\PBP(H)$ has an extended formulation with at most $(\frac{2}{3} 3^{r}+2(r-1)(2^r+1))|V|$ variables and inequalities. In particular, if $k= O(\log \poly(|V|, |S|))$, then $\PBP(H)$ has a polynomial-size extended formulation. Moreover, all coefficients and right-hand side constants in the system defining $\PBP(H)$ are $0,\pm 1$.
\end{theorem}

\begin{prf}
By definition of $\alpha$-acyclic hypergraphs, the underlying hypergraph of $H$ has a sequence of $\alpha$-leaves of length $n:=|V|$. Let us denote this sequence by $v_1,v_2,\dots,v_n$.
Since $v_1$ is an $\alpha$-leaf of the underlying hypergraph of $H$, there exists a maximal (for inclusion) signed edge $\bar s$ of $H$ containing $v_1$.
Let $H'$ be obtained from $H$ by first inflating each signed edge $s'$ containing $v_1$ such that $s' \subset \bar s$ to $\bar s$, and then by adding all signed edges $s - v_1$ for all signed edges containing $v_1$. 
Denote by $N_{v_1}$ the number of signed edges $s \in S \setminus \{\bar s\}$ containing $v_1$ such that $s \subset \bar s$.
By~\cref{lem inflation}, it suffices to find an extended formulation of $\PBP(H')$. An extended formulation for $\PBP(H)$ is then obtained by juxtaposing the extended formulation for $\PBP(H')$ and $N_{v_1}$ equalities 
containing at most $N_{v_1}$ additional variables.
Since $H$ is a rank-$r$ signed hypergraph we have:
$$
N_{v_1}\leq N_{\rm max} = 2\Big(\sum_{k=1}^{r-2}\binom{r-1}{k} 2^{r-1-k}\Big)=2(3^{r-1}-2^{r-1}-1).
$$
%
Moreover, by construction, $v_1$ is a $\beta$-leaf of of the underlying hypergraph of $H'$, 
implying we can apply the decomposition result of~\cref{th decomp}.
Namely, let $V_1$ be the underlying edge of $\bar s$, let $S_{v_1}$ be the set of signed edges of $H'$ parallel to $\bar s$, and let 
$P_{v_1} := \{s-v_1 : |s-v_1| \ge 2\}$.
Define the signed hypergraphs $H_{v_1}:=(V_1, S_{v_1} \cup P_{v_1})$ and $H'_1 = H'-v_1$.
Then, by \cref{th decomp}, $\PBP(H')$ is decomposable into $\PBP(H_{v_1})$ and $\PBP(H'_1)$.
An extended formulation for the pseudo-Boolean polytope of $H_{v_1}$ can then be obtained from Part~2 of Remark~\ref{rem1}.

Next, we show that we can apply the above construction recursively to $H' - v_1$ with sequence of $\alpha$-leaves $v_2,\dots,v_n$.
Note that the underlying hypergraph of $H' - v_1$ is obtained from the underlying hypergraph of $H - v_1$ by possibly removing some edges contained in $\bar s - v_1$.
First, note that the rank of $H' - v_1$ is at most $r$, as removing nodes and edges from a hypergraph cannot increase its rank.
Next, we show that $v_2,\dots,v_n$ is a sequence of $\alpha$-leaves of the underlying hypergraph of $H' - v_1$.
For $k \in \{1,\dots,n\}$, we show that $v_k$ is an $\alpha$-leaf of the underlying hypergraph of $H'_k := H' - v_1 - v_2 - \cdots - v_{k-1}$.
Since $v_k$ is an $\alpha$-leaf of the underlying hypergraph of $H_k := H - v_1 - v_2 - \cdots - v_{k-1}$, let $\hat s$ be a signed edge of $H_k$ containing $v_k$ that is maximal for inclusion.
Let $\tilde s$ be a signed edge of $H$ such that $\hat s = \tilde s - v_1 - v_2 - \cdots - v_{k-1}$.
If $\tilde s$ is also a signed edge of $H'$ we are done.
Otherwise, by definition of $H'$, we have $\tilde s \subseteq \bar s$.
But then $\bar s' := \bar s - v_1 - v_2 - \cdots - v_{k-1}$ is a signed edge of $H'_k$ and we have $\hat s \subseteq \bar s'$.
Therefore, $\bar s'$ is a signed edge of $H'_k$ containing $v_k$ that is maximal for inclusion.
Hence, we can apply the above construction recursively.
The proof follows by induction on $n$.

Define the signed hypergraph $H_{v_k}=(V_k, S_{v_k}\cup P_{v_k})$, $k \in \{1,\dots,n\}$,
where $S_{v_k}$ is a set of parallel signed edges containing $v_k$, and $P_{v_k}:=\{s-v_k: |s-v_k|\geq 2\}$. Moreover, $V_k$ is the underlying edge of a signed edge in $S_{v_k}$. By the above argument an
extended formulation for $\PBP(H)$ is obtained by juxtaposing extended formulations of $\PBP(H_{v_k})$ for all $k \in \{1,\dots,n\}$ and at most $N_{\max}|V|$ equalities containing at most $N_{\max}|V|$ additional variables. 

Finally, consider the pointed signed hypergraph $H_{v_k}$ for some $k \in \{1,\dots,n\}$. Since the rank of $H$ is $r$, we have $|S_{v_k}| \leq 2^r$. Moreover, since by construction all signed edges in $S_{v_k}$ are parallel, by Part~(ii) of Remark~\ref{rem3} we conclude that
$\PBP(H_{v_k})$ has an extended formulation with at most $(r-1) 2^r+ r$ variables and at most 
$2((r-1) 2^r+r)$ inequalities.
%
Therefore, the pseudo-Boolean polytope $\PBP(H)$ has an extended formulation with at most
$$
((r-1) 2^r+ r)|V| + 2(3^{r-1}-2^{r-1}-1)|V|=
(\frac{2}{3} 3^{r}+(r-2) (2^{r}+1))|V|,
$$
variables, and at most
$$
(2(r-1) 2^r+2 r)|V|+ 2(3^{r-1}-2^{r-1}-1)|V| \leq 
(\frac{2}{3} 3^{r}+2(r-1)(2^r+1))|V|
$$
inequalities.
\end{prf}

\begin{remark}
The proof of Theorem~\ref{alphaStuff general} follows from the RID framework. To see this, first note that the node $\bar v$ chosen in Step~2 of the $i$th iteration of the RID framework is the $\alpha$-leaf $v_i$ for all $i \in [n]$. 
Since by definition of $\alpha$-leaves each $v_i$ is contained in a maximal edge $\bar s$, the inflation operation of Step 2 is defined as follows: inflate all edges containing $v_i$ to $\bar s$. Since $H$ is a rank-$r$ hypergraph, from the proof of Theorem~\ref{alphaStuff general} it follows that the total number of edges added due to inflation operations is upper bounded by $2^r|V|$, which is a polynomial in $|V|,|S|$ if we have $r= O(\log\poly(|V|,|S|))$.
\end{remark}

\subsection{Hypergraphs with log-poly gaps}
\label{sec: smallgap}

Consider a hypergraph $G=(V,E)$, and let $E' \subseteq E$; we define the \emph{gap} of $E'$ as
$$
\gap(E') := \max\Bigg\{ \Big|\bigcup_{f \in E'} f \Big| - |e| : e \in E' \Bigg\}.
$$
Moreover, for a signed hypergraph $H=(V,S)$, and $S' \subseteq S$, the gap of $S'$, denoted by $\gap(S')$ is defined as the gap of the set of the underlying edges of $s \in S'$. 

Consider the signed hypergraph $H=(V,S)$ defined in the statement of Proposition~\ref{cor1}. It can be checked that the gap of $S$ is upper bounded by $k$. Then, from the proof of Proposition~\ref{cor1} it follows that, if $\gap(S) = O(\log\poly(|V|,|S|))$, the pseudo-Boolean polytope $\PBP(H)$ has a polynomial-size extended formulation.
%
In this section, we strengthen this result by providing polynomial-size extended formulations for the pseudo-Boolean polytope of signed hypergraphs $H=(V,S)$
with $\gap(S')= O(\log\poly(|V|,||S|))$
for only certain (small) subsets $S' \subset S$.


In the following, given a signed hypergraph $H=(V,S)$, we say that $s \in S$ is a \emph{maximal signed edge} if there is no $s' \in S$ with $s' \supset s$. 

\begin{proposition}\label{cor2}
Consider a signed hypergraph $H=(V,S)$ of rank $r$. For each $s \in S$, among all maximal signed edges of $H$ containing $s$, denote by $f_s$ one with minimum cardinality. Let $k$ be such that 
\begin{equation}\label{gcond}
\gap(\{s,f_s\}) \leq k, \quad \forall s \in S.    
\end{equation}
Denote by $\bar S$ the set of maximal signed edges of $H$.
If the underlying hypergraph of $(V, \bar S)$ is $\beta$-acyclic, then $\PBP(H)$ has an extended formulation with $O(r 2^k|V||S|)$ variables and inequalities.
In particular, if $k= O(\log \poly(|V|, |S|))$, then $\PBP(H)$ has a polynomial-size extended formulation. Moreover, all coefficients and right-hand side constants in the system defining $\PBP(H)$ are $0,\pm 1$.
\end{proposition}

\begin{prf}
Denote by $H'=(V,S')$ the signed hypergraph obtained from $H$ by inflating every non-maximal signed edge $s \in S$ to $f_s$, where $f_s$ is defined in the statement. From~\eqref{gcond} it follows that $|S'|\leq 2^{k}(|S|-|\bar S|)+|\bar S| \leq 2^k |S|$.
Notice that the underlying hypergraph of $H'$ coincides with the underlying hypergraph of $(V, \bar S)$, which by assumption is $\beta$-acyclic. Hence by Theorem~\ref{extended}, $\PBP(H')$ has an extended formulation with at most $2r|V|(2^{k+1} |S|+1)$ variables, and at most $4r|V|(2^{k+1} |S|+1)$ inequalities. Therefore, from Lemma~\ref{lem inflation}
we deduce that an extended formulation for $\PBP(H)$ is obtained by adding at most $|S|-|\bar S|$ equalities consisting of at most $|S|-|\bar S|$ additional variables to $\PBP(H')$, and this completes the proof.
\end{prf}

The next example illustrates the power of~\cref{cor2}.

\begin{example}\label{expl1}
    Consider the hypergraph $G=(V,E)$ with $V=\{v_1, \ldots, v_{2n}\}$ for some $n \geq 3$ and $E=\{e_1, \ldots, e_{n}\} \cup \{f_1, \ldots, f_{n}\} \cup \{g_1, \ldots, g_{n}\}$, where $e_1=\{v_1, v_2, v_3\}$, $e_i=\{v_{2i-2}, v_{2i-1}, v_{2i}, v_{2i+1}\}$ for all $i \in \{2,\ldots,n-1\}$, $e_n =\{v_{2n-2}, v_{2n-1}, v_{2n}\}$, $f_1=\{v_1, v_2\}$, $f_i = \{v_{2i-2},v_{2i}\}$ for all $i \in \{2,\ldots,n\}$, $g_i =\{v_{2i-1}, v_{2i+1}\}$ for all $i \in \{1, \ldots, n-1\}$ and $g_n=\{v_{2n-1}, v_{2n}\}$. See Figure~\ref{fig long} for an illustration. It can be checked that $G$ has a $\beta$-cycle of length $2n$ given by 
    $$C= v_1, f_1, v_2, f_2, v_4, \ldots ,v_{2n-2}, f_n, v_{2n}, g_n, v_{2n-1}, \ldots ,v_5, g_2, v_3, g_1, v_1.$$ We now use Proposition~\ref{cor2} to obtain a polynomial-size extended formulation for $\MP(G)$. To this end, notice that $\gap(\{f_1, e_1\}) = \gap(\{g_1, e_1\}) = \gap(\{f_n, e_n\}) = \gap(\{g_n, e_n\}) = 1$ and $\gap(\{f_i, e_i\}) = \gap(\{g_i, e_i\}) = 2$ for all $i \in \{2,\ldots,n-1\}$. Moreover hypergraph $(V, \bar E)$, where $\bar E = \{e_1, \ldots, e_n\}$ is $\beta$-acyclic. Therefore, by proof of Proposition~\ref{cor2}, by inflating each $f_i, g_i$ to $e_i$ for all $i \in \{1,\ldots,n\}$, one obtains a signed hypergraph $H=(V,S)$ whose underlying hypergraph  is $\beta$-acyclic and $|S| \leq 4|E|$. Therefore, by Theorem~\ref{extended}, $\PBP(H)$ and as a result $\MP(G)$ have polynomial-size extended formulations.

    We should remark that while $G$ is not $\beta$-acyclic, it is an $\alpha$-acyclic hypergraph with rank $r=4$. Hence, Theorem~\ref{alphaStuff general} also provides
    (a different) polynomial-size extended formulation for $\MP(G)$. However, it is simple to modify this example so that $G$ will no longer have a log-poly rank and hence 
    Theorem~\ref{alphaStuff general} will no longer be applicable. Denote by $U_i$ a set of $n$ nodes such that $U_i \cap V = \emptyset$ for all $i \in [n-1]$
    and $U_i \cap U_j = \emptyset$ for all $i \neq j \in [n-1]$.
    Define the hypergraph $G'=(V',E')$, with $V'=\bigcup_{i\in [n-1]}{U_i} \cup V$ and $E'=\{e'_1, \ldots, e'_n\} \cup \{f'_1, \ldots, f'_n\} \cup \{g'_1, \ldots, g'_n\}$, where $e'_1= e_1 \cup U_1$, $f'_1= f_1 \cup U_1$, $g'_1= g_1 \cup U_1$, 
    $e'_i = e_i \cup U_{i-1} \cup U_{i}$,
    $f'_i = f_i \cup U_{i-1} \cup U_{i}$, 
    $g'_i = g_i \cup U_{i-1} \cup U_{i}$ for $i \in \{2,\ldots,n-1\}$ and $e'_n= e_n \cup U_{n-1}$, $f'_n= f_n \cup U_{n-1}$, $g'_n= g_n \cup U_{n-1}$.
    Then the rank of $G'$ is larger than $2 n$ and hence 
    Theorem~\ref{alphaStuff general} is no longer applicable. However, the gaps remain unchanged; \ie $\gap(\{f_i, e_i\}) = \gap(\{f'_i, e'_i\})$ and $\gap(\{g_i, e_i\}) = \gap(\{g'_i, e'_i\})$ for all $i \in [n]$ and the hypergraph $(V', \bar E)$ is $\beta$-acyclic, implying that by Proposition~\ref{cor2} one can obtain a polynomial-size extended formulation for $\MP(G')$. 
\end{example}

\begin{figure}[h]
\begin{center}
\includegraphics[width=.9\textwidth]{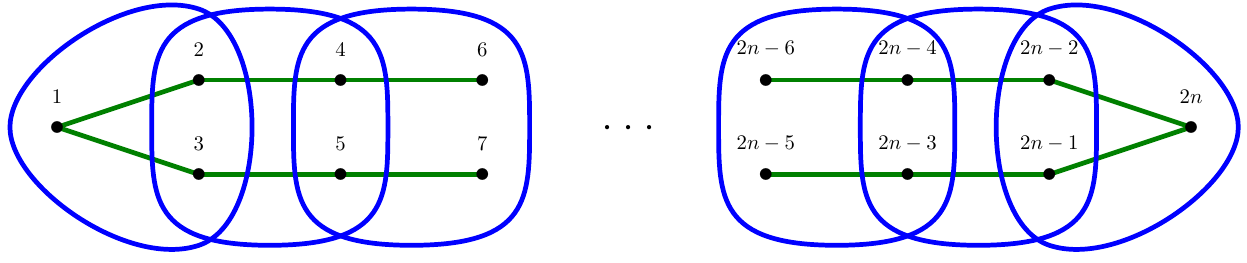}
\caption{The hypergraph $G=(V,E)$ of Example~\ref{expl1}; $G$ contains a long $\beta$-cycle of length $|V|$.}
\label{fig long}
\end{center}
\end{figure}

Let $G=(V, E)$ be a hypergraph and let 
$C = v_1, e_1, v_2, e_2, \dots , v_q, e_q, v_1$ for some $q \geq 3$ be a $\beta$-cycle  
of  $G$. Then the \emph{support hypergraph} of $C$ is the hypergraph $G[C]=(V[C], E[C])$, where $V[C] = \bigcup_{i=1}^q e_i$ and $E[C] = \{e_1,\dots,e_q\}$. The next proposition essentially indicates that if edge sets of the support hypergraphs of $\beta$-cycles of the underlying hypergraph of $H$ have log-poly gaps, then one can \emph{remove} these $\beta$-cycles by inflating the edges of the cycle, and obtain a polynomial-size extended formulation for $\PBP(H)$.

\begin{proposition}
\label{cor3}
Consider a signed hypergraph $H=(V,S)$ of rank $r$, and denote by $G = (V,E)$ its underlying hypergraph.
Let $\cal C$ denote the set of all $\beta$-cycles in $G$, and for each $C \in \cal C$, let $G[C]=(V[C], E[C])$ be the support hypergraph of $C$.
Let $\tilde G$ be the hypergraph $(\bigcup_{C \in \cal C} V[C], \bigcup_{C \in \cal C} E[C])$,
let $(V_1,E_1), (V_2,E_2),\dots,(V_\omega,E_\omega)$ for some $\omega \geq 1$ be the connected components of $\tilde G$, and 
let $k$ be such that
\begin{equation}\label{gcond2}
\gap(E_i) \leq k, \quad \forall  i \in [\omega].
\end{equation}
Then, $\PBP(H)$ has an extended formulation with $O(r 2^k |V| |S|)$ variables and  inequalities.
In particular, if $k=O(\log \poly(|V|, |S|))$, then $\PBP(H)$ has a polynomial-size extended formulation. Moreover, all coefficients and right-hand side constants in the system defining $\PBP(H)$ are $0,\pm 1$.
\end{proposition}

\begin{prf}
Each $\beta$-cycle in $G$ is contained in precisely one of the connected components of $\tilde G$.
Let $H'$ be obtained from $H$ by inflating, for each $j \in \{1,2,\dots,\omega\}$, every signed edge $s$ contained in $V_j$, to $V_j$.
Let $G'$ be the underlying hypergraph of $H'$.
We show that $G'$ contains no $\beta$-cycle.
Assume for a contradiction that $G'$ contains a $\beta$-cycle $C'$.
Consider first the case where $C'$ does not contain any set $V_j$ as an edge.
Then $C'$ is a $\beta$-cycle in $G$ not contained in any connected component of $\tilde G$, which gives us a contradiction.
Next, consider the case where $C'$ contains at least one set $V_j$ as an edge.
We now show how we can modify $C'$ to obtain a $\beta$-cycle $C$ in $G$, not contained in any connected component of $\tilde G$, which gives us again a contradiction.
Assume that $C'$ contains the edge $V_j$, and let $u,v$ be the nodes in $C'$ before and after $V_j$.
We then define $C$ by replacing, in $C'$, the edge $V_j$ with a minimal sequence $f_1,w_2,\dots,w_t,f_t$, where $f_j$ are edges in $E_j$ and $w_j$ are nodes in $V_j$, such that $u \in f_1$ and $v \in f_t$. 
Note that this sequence exists because $(V_j,E_j)$ is connected.
Clearly, nodes $w_2,\dots,w_t$ and edges $f_1,\dots,f_t$ did not belong to $C'$.
It is simple to check that applying recursively the above construction to all sets among $V_1,\dots,V_\omega$ contained in $C'$ yields a $\beta$-cycle in $C$ not contained in any connected component of $\tilde G$, which gives us a contradiction.
We have therefore shown that $G'$ is $\beta$-acyclic.

By Lemma~\ref{lem inflation} and assumption~\eqref{gcond2}, an extended formulation for $\PBP(H)$ can be obtained by juxtaposing an extended formulation for
$\PBP(H')$ together with at most 
$\bigcup_{C \in \cal C} |E[C]| \leq |S|$ 
equalities, containing at most 
$\bigcup_{C \in \cal C} |E[C]| \leq |S|$
additional variables.
Now consider the signed hypergraph $H'=(V,S')$.
We then have $|S'| \leq |S|-\sum_{i \in [\omega]}{|E_i|}+\sum_{i \in [\omega]}{2^{k}|E_i|} \leq 2^k |S|$.
Since the underlying hypergraph of $H'$ is $\beta$-acyclic, by Corollary~\ref{extended},  $\PBP(H')$ has an extended formulation with at most $2r|V|(2^{k+1} |S|+1)$ variables, and at most $4r|V|(2^{k+1} |S|+1)$ inequalities.
\end{prf}




\begin{remark}\label{becor3}
Under more restrictive assumptions, the removal of $\beta$-cycles can be performed in a more efficient manner than the technique in the proof of \cref{cor3}. 
Two $\beta$-cycles are \emph{equivalent} if one can be obtained from the other via a circular permutation and/or reversing the order of the edges.
Now consider the case where all $\beta$-cycles in each $(V_i, E_i)$, $i \in [\omega]$, are equivalent. Also notice that each $\beta$-cycle is contained in precisely one of the connected components
$(V_i, E_i)$.
For each $i \in [\omega]$, denote by $f_i$ an edge of maximum cardinality in $E_i$.  By the definition of a $\beta$-cycle, each node of the cycle should be present in exactly two edges of the cycle. Therefore, by inflating $f_i$ to $V_i$, we can eliminate all $\beta$-cycles in this connected component. 
Let $k$ be such that
\begin{equation}\label{wass}
\gap(\{f_i, V_i\}) \leq k, \quad \forall C \in \cal C.
\end{equation}
Then one can replace assumption~\eqref{gcond2} by the weaker assumption~\eqref{wass}
\end{remark}

Let 
$C = v_1, e_1, v_2, e_2, \dots , v_q, e_q, v_1$ be a $\beta$-cycle of length $q$ for some $q \geq 3$ in $G$. Define $\delta = \max_{i\in \{1,\ldots,q\}}{|e_i|}$. We then have $\gap(E[C]) \geq q- \delta$.
Therefore, Proposition~\ref{cor3} (or its enhancement described in Remark~\ref{becor3}) is not applicable for hypergraphs containing long $\beta$-cycles consisting of edges of small cardinality.  For example, the hypergraph $G$ of Example~\ref{expl1} does not satisfy the assumptions of Proposition~\ref{cor3}, since the gap of the edge set of its $\beta$-cycle equals $2n-2$, which in turn equals $|V|-2$. The following example illustrates the usefulness of Proposition~\ref{cor3}.  

\begin{example}\label{expl2}
Consider a hypergraph $G=(V,E)$ with $V=\{v_1,\ldots, v_n\}$ for some $n \geq 12$ and 
\begin{align*}
E=\Big\{V, & \{v_1, v_2\}, \{v_2, v_3\}, \{v_4, v_5\},
\{v_5, v_6\}, \{v_4, v_7\}, \{v_6, v_7\}, \{v_4,v_5, v_6\}, \{v_8, v_9\}, \{v_9, v_{10}\}, \{v_8, v_{10}\}, \\
& \{v_9,v_{11}, v_{12}\}, \{v_{10},v_{11}, v_{12}\}\Big\}.
\end{align*}
See Figure~\ref{fig inflate_prop4} for an illustration. It then follows that the hypergraph $\tilde G$ defined in Proposition ~\ref{cor3} 
has two connected components: $(V_1, E_1)$
where $V_1 = \{v_4, v_5, v_6, v_7\}$, $E_1=\{e \in E: e \subseteq V_1\}$, and $(V_2, E_2)$ where $V_2 =\{v_8, v_9, v_{10}, v_{11}, v_{12}\}$, $E_2=\{e \in E: e \subseteq V_2\}$. We then have $\gap(E_1) = 2$ and $\gap(E_2)=3$. Hence the assumptions of Proposition~\ref{cor3} are satisfied and by inflating every edge in $E_1$ to $V_1$ and by inflating every edge in $E_2$ to $V_2$, we obtain a signed hypergraph $H$ whose underlying hypergraph is $\beta$-acyclic. Therefore, by Theorem~\ref{extended}, $\PBP(H)$ and as a result $\MP(G)$ have polynomial-size extended formulations.
\end{example}

Notice that the hypergraph of Example~\ref{expl2}
does not satisfy the assumptions of Proposition~\ref{cor2} as we have $\gap(\{s, f_s\})=|V|-2$ for $s = \{v_1, v_2\}$. That is, Propositions~\ref{cor2} and~\ref{cor3} are rather \emph{complementary} and an effective combination of the two can be used to obtain polynomial-size extended formulations for many signed hypergraphs whose underlying hypergraphs contain $\beta$-cycles.

\begin{figure}[h]
\begin{center}
\includegraphics[width=.7\textwidth]{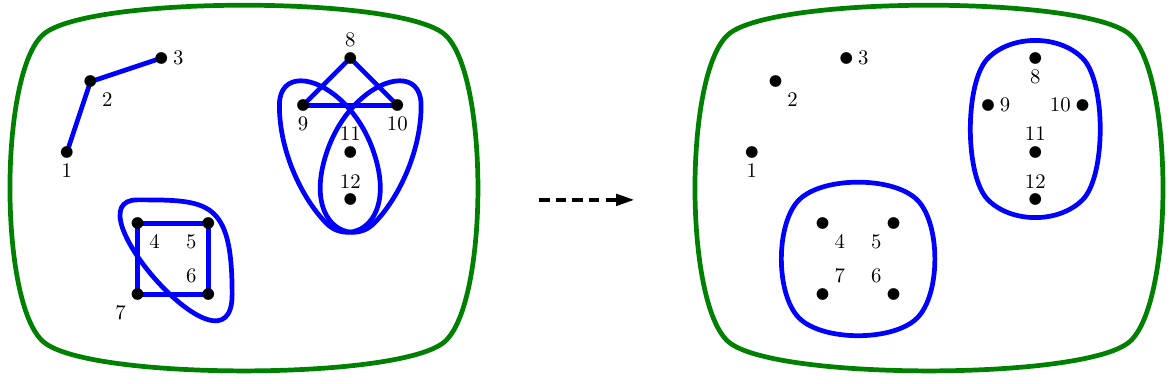}
\caption{The hypergraph $G$ of Example~\ref{expl2} with $n=12$ is shown on the left. The $\beta$-cyclic hypergraph obtained from $G$ after inflation operations is shown on the right.}
\label{fig inflate_prop4}
\end{center}
\end{figure}

\begin{remark}
Proofs of Propositions~\ref{cor2} and~\ref{cor3} follow from Step 1 of the RID framework: given a signed hypergraph $H=(V,S)$, via a number of inflation operations, we obtain the hypergraph $H'=(V,S')$ such that: 
\begin{itemize}
\item [(i)] the underlying hypergraph of $H'$ is $\beta$-acyclic, hence, by Theorem~\ref{extended}, the polytope $\PBP(H')$ has a polynomial-size extended formulation in $|V|, |S'|$, 

\item [(ii)] $|S'| \leq |S| + \poly(|V|,|S|)$, where this inequality is satisfied because of the log-poly gaps.
\end{itemize}
\end{remark}

\subsection{A more general framework}
\label{sec: general}

So far, all inflation operations we employed to obtain polynomial-size extended formulations  (\ie Theorem~\ref{alphaStuff general}, and Propositions~\ref{cor1}--\ref{cor3}) have a simple form: we inflated a number of signed edges to the same set. This in turn enabled us to use Theorems~\ref{th decomp} and ~\ref{pointedHull} to obtain our extended formulations. As we detail next, the inflation of signed edges can be useful in a more general setting. Consider a signed hypergraph $H=(V,S)$. Assume that $S$ is not totally ordered. Let $\E$ denote a set of subsets of $V$ that is totally ordered. Moreover suppose that for each $s \in S$, there exists some $e \in \E$ such that $e \supseteq s$. Now inflate each $s \in S$ to some element in $\E$, which we denote by $e(s)$.
 Define 
 $$I(S,\E):=\bigcup_{s \in S}{I(s, e(s))}.$$ 
 Define the signed hypergraph $H'=(V,I(S,\E))$. By Part~2
 of Remark~\ref{rem1}, the pseudo-Boolean polytope $\PBP(H')$ has a polynomial-size extended formulation in $|V|, |I(S,\E)|$. Therefore, from Lemma~\ref{lem inflation} we deduce that if $|I(S,\E)|$ is a polynomial in $|V|,|S|$, we can obtain a polynomial-size extended formulation for $\PBP(H)$
 as well. 
 The following proposition is an illustration of how this general setting can be used to obtain a polynomial-size extended formulation for the pseudo-Boolean polytope.

\begin{proposition}\label{cor4}
Consider a signed hypergraph $H=(V,S)$. Let $S=S_1 \cup S_2$, such that each $s \in S_1$ contains at least $|V|-k_1$ nodes, whereas each $s \in S_2$ is contained in $U \subset V$, where $|U|=k_2$. Then $\PBP(H)$ has an extended formulation with $O(2^{\max\{k_1,k_2\}}|S||V|)$ variables and inequalities.
In particular, if $\max\{k_1,k_2\}= O(\log \poly(|V|, |S|))$, then $\PBP(H)$ has a polynomial-size extended formulation. Moreover, all coefficients and right-hand side constants in the system defining $\PBP(H)$ are $0,\pm 1$.
\end{proposition}

\begin{prf} 
Consider the signed hypergraph $H'=(V,S')$ obtained from $H$ by inflating every $s \in S_1$ with $|s| < |V|$ to $V$, and by inflating every $s \in S_2$ with $|s| < |U|$ to $U$. We then have $|S'|\leq 2^{k_1}|S_1|+2^{k_2}|S_2| \leq 2^{\max\{k_1,k_2\}}|S|$. The remainder of the proof is identical to the proof of Proposition~\ref{cor1}.
\end{prf}

We conclude the paper by remarking that while we presented several classes of signed hypergraphs for which our proposed RID framework provides polynomial-size extended formulations for the corresponding pseudo-Boolean polytopes, a complete characterization of such signed hypergraphs remains an open question and is a subject of future research.

\section{Appendix}
\label{sec: appendix}

\paragraph{Proof of Proposition~\ref{nestedHull}.}
The proof is by induction on the number of nodes $n$ of $H$. The base case is $n=2$; in this case we have $S = \E_{1} = \{q_1, q_2, q_3, q_4\}$. It is simple to check that $\PBP(H)$ is given by:
\begin{align*}
   & z_{q_1} + z_{q_2} = z_{v_1}, \; 
    z_{q_1} + z_{q_3} = z_{v_2}, \; 
    z_{q_3} + z_{q_4} = 1-z_{v_1}\\
   & z_{q_1}, z_{q_2}, z_{q_3}, z_{q_4} \geq 0.
\end{align*}
Henceforth, let $n \geq 3$. Denote by $\PBP^0(H)$ 
(resp. $\PBP^1(H)$) the face of $\PBP(H)$ defined by $z_{v_{n}} = 0$ (resp. $z_{v_n} = 1$). 
We then have: 
\begin{equation}\label{disjp}
\PBP(H) = \conv\Big(\PBP^0(H) \cup \PBP^1(H)\Big).    
\end{equation}
Define the signed hypergraph $\bar H := (\bar V, \bar S)$ with $\bar V := V \setminus \{v_n\}$ and $\bar S := S \setminus \E_{n-1}$. Since $\bar H$ has one fewer node than $H$ and is a nested signed hypergraph, by the induction hypothesis, the pseudo-Boolean polytope $\PBP(\bar H)$ is given by:
   \begin{align*}
   & z_{s} + z_{\ell(s)} = z_{p(s)}, \quad \forall s \in \E^{+}_k, k \in [n-2]\setminus\{1\}\nonumber\\
   & z_s \geq 0, \quad \forall s \in S \setminus \E_{n-1}\nonumber\\  
    &     \sum_{s \in \E^{+}_k}{z_s}\leq z_{v_{k+1}}, \quad \forall k \in [n-2] \setminus \{1\}\\
     &    \sum_{s \in \E^{-}_k}{z_s}\leq 1-z_{v_{k+1}}, \quad \forall k \in [n-2] \setminus \{1\}\nonumber\\
      &      z_{q_1} + z_{q_2} = z_{v_1}, \; 
    z_{q_1} + z_{q_3} = z_{v_2}, \; 
    z_{q_3} + z_{q_4} = 1-z_{v_1}. \nonumber
    \end{align*}
Denote by $\bar z$ the vector consisting of $z_v$ for all $v \in \bar V$, and $z_s$ for all $s \in \bar S$.
It then follows that
\begin{align*}
& \PBP^0(H) = \Big\{z \in \R^{V \cup S}: z_{v_n} = 0, \; z_s =0,  \; \forall s \in \E^{+}_{n-1}, z_s = z_{p(s)}, \; \forall e \in \E^{-}_{n-1},\; \bar z \in \PBP(\bar H) \Big\}\\
& \PBP^1(H) = \Big\{z \in \R^{V \cup S}: z_{v_n} = 1, \; z_s =z_{p(s)},  \; \forall e \in \E^{+}_{n-1}, z_s = 0, \; \forall s \in \E^{-}_{n-1}, \; \bar z \in \PBP(\bar H) \Big\}.
\end{align*}
By~\eqref{disjp} and using Balas' formulation for the union of polytopes~\cite{Bal98}, it follows that $\PBP(H)$ is the projection onto the space of the $z_v$, $v \in V$, and $z_s$, $s \in S$, variables of the polyhedron defined by the following system~\eqref{disjunc1}--\eqref{disjunc3}:
\begin{align}
\begin{split}\label{disjunc1}
   & \lambda_0 + \lambda_1 = 1, \; \lambda_0 \geq 0\;, \lambda_1 \geq 0\\
   &  z_v = z^0_v + z^1_v, \quad \forall v \in V\\
   &  z_s = z^0_s + z^1_s, \quad \forall s \in S 
\end{split}
\end{align}

\begin{align}
\begin{split}\label{disjunc2}
   & z^0_{v_n} = 0\\
   & z^0_s = 0, \quad \forall s \in \E^{+}_{n-1}\\
   &  z^0_{s} = z^0_{p(s)}, \quad \forall s \in  \E^{-}_{n-1}\\
   &  z^0_{s} + z^0_{\ell(s)} = z^0_{p(s)},  \quad \forall s \in \E^{+}_k, k \in [n-2]\setminus\{1\}\\
   & z^0_s \geq 0,  \quad \forall s \in S \setminus \E_{n-1}\\
   &  \sum_{s \in \E^{+}_k}{z^0_s}\leq z^0_{v_{k+1}}, \quad \forall k \in [n-2] \setminus \{1\}\\
   &   \sum_{s \in \E^{-}_k}{z^0_s}\leq \lambda_0-z^0_{v_{k+1}}, \quad \forall k \in [n-2] \setminus \{1\}\\
   &  z^0_{q_1} + z^0_{q_2} = z^0_{v_1}, \; 
    z^0_{q_1} + z^0_{q_3} = z^0_{v_2}, \; 
    z^0_{q_3} + z^0_{q_4} = \lambda_0-z^0_{v_1}
\end{split}
\end{align}

\begin{align}
\begin{split}\label{disjunc3}
   & z^1_{v_n} = \lambda_1\\
   & z^1_s = z^1_{p(s)}, \quad \forall s \in \E^{+}_{n-1}\\
   &  z^1_{s} = 0, \quad \forall s \in \E^{-}_{n-1}\\
   &     z^1_{s} + z^1_{\ell(s)} = z^1_{p(s)}, \quad \forall s \in \E^{+}_k, k \in [n-2]\setminus\{1\}\\
   &  z^1_s \geq 0,\quad \forall s \in S \setminus \E_{n-1}\\
   &  \sum_{s \in \E^{+}_k}{z^1_s} \leq z^1_{v_{k+1}}, \quad \forall k \in [n-2] \setminus \{1\}\\
   &  \sum_{s \in \E^{-}_k}{z^1_s} \leq \lambda_1-z^1_{v_{k+1}}, \quad \forall k \in [n-2] \setminus \{1\}\\
   & z^1_{q_1} + z^1_{q_2} = z^1_{v_1}, \; 
    z^1_{q_1} + z^1_{q_3} = z^1_{v_2}, \; 
    z^1_{q_3} + z^1_{q_4} = \lambda_1-z^1_{v_1}.
    \end{split}
\end{align}
In the remainder of this proof, we project out $z^0, z^1, \lambda_0, \lambda_1$ from system~\eqref{disjunc1}--\eqref{disjunc3} and obtain a description of $\PBP(H)$ in the original space.
From $z_{v_n} = z^0_{v_n}+ z^1_{v_n}$, $z^0_{v_n} = 0$,
$z^1_{v_n} = \lambda_1$, and $\lambda_0+\lambda_1=1$, it follows that
\begin{equation}\label{lss}
\lambda_0 = 1- z_{v_n}, \quad \lambda_1 = z_{v_n}.    
\end{equation}
By $z_s = z^0_s + z^1_s$ for all $s \in S$, $z^0_s =0$, and $z^1_s = z^1_{p(s)}$  for all $s \in \E^{+}_{n-1}$, we get 
\begin{equation}\label{proj1}
z^1_s = z^1_{p(s)} = z_s, \quad z^0_{p(s)}=z_{p(s)} -z_s, \quad \forall s \in \E^{+}_{n-1}.
\end{equation}
Similarly, by
$z^1_{\ell(s)} =0$ and $z^0_{\ell(s)} = z^0_{p(s)}$  for all $s \in \E^{+}_{n-1}$, we get 
\begin{equation}\label{proj2}
z^0_{\ell(s)} = z^0_{p(s)}=z_{\ell(s)}, \quad  z^1_{p(s)}=z_{p(s)} -z_{\ell(s)}, \quad \forall s \in \E^{+}_{n-1}.
\end{equation}
Using~\eqref{disjunc1} to project out $z^0_v$, $v \in V$, $z^0_s$, $s \in S$, using~\eqref{lss} to project out $\lambda_0, \lambda_1$, and using~\eqref{proj1} and~\eqref{proj2} to project out $z^0_s,z^1_s$, $s \in \E_{n-1}$ with $N(s) \neq \emptyset$, we deduce that system~\eqref{disjunc1}-\eqref{disjunc3} simplifies to:
\begin{align}
    & z_{s + v^-_n} +  z_{\ell(s) + v^-_n} = z_{p(s)}-z^1_{p(s)}, \; \forall s \in \E^{+}_{n-2}: s + v^-_n, \ell(s) + v^-_n \in \E^{-}_{n-1}\label{n11}\\
    & z_{s + v^-_n} +  z_{\ell(s)}-z^1_{\ell(s)} = z_{p(s)}-z^1_{p(s)}, \; \forall s \in \E_{n-2}: s + v_n^- \in \E^{-}_{n-1},  N(\ell(s))=\emptyset\label{n12}\\
    & z_s-z^1_{s} + z_{\ell(s)}-z^1_{\ell(s)} = z_{p(s)}-z^1_{p(s)}, \; \forall s \in \E^{+}_k: k \in [n-3]\setminus\{1\}, \; {\rm or} \nonumber\\
    & \qquad \qquad \qquad \qquad \qquad \qquad \qquad \qquad \qquad \qquad k=n-2, N(s)=N(\ell(s)) = \emptyset\label{n14}\\
    & z_s - z^1_s \geq 0,\quad \forall s \in \E_{k}: k \in [n-3] \; {\rm or}  \; k = n-2, N(s) = \emptyset\label{n15}\\
    & \sum_{\substack{s \in \E^{+}_{n-2}: \\ s + v^{-}_n \in \E^{-}_{n-1}}}{z_{s+ v^-_n}}+
    \sum_{\substack{s \in \E^{+}_{n-2}: \\ N(s) = \emptyset}}{(z_s-z^1_s)}
    \leq z_{v_{n-1}}-z^1_{v_{n-1}}, \label{n16}\\
     & \sum_{\substack{s \in \E^{-}_{n-2}: \\ s+ v^-_n \in \E^{-}_{n-1}}}{z_{s+v^-_n}}+
     \sum_{\substack{s \in \E^{-}_{n-2}: \\ N(s) = \emptyset}}{(z_s-z^1_s)}\leq 1-z_{v_n}-z_{v_{n-1}}+z^1_{v_{n-1}}, \label{n17}\\
      & \sum_{s \in \E^{+}_k}{(z_s-z^1_s)}\leq z_{v_{k+1}}-z^1_{v_{k+1}}, \quad \forall k \in [n-3] \setminus \{1\}\label{n18}\\
     & \sum_{s \in \E^{-}_k}{(z_s-z^1_s)}\leq 1-z_{v_n}-z_{v_{k+1}}+z^1_{v_{k+1}}, \quad \forall k \in [n-3] \setminus \{1\}\label{n19}\\
     \begin{split}\label{n19p}
    & z_{q_1}-z^1_{q_1} + z_{q_2}- z^1_{q_2} = z_{v_1}-z^1_{v_1}, \; 
    z_{q_1}-z^1_{q_1} + z_{q_3}- z^1_{q_3} = z_{v_2}-z^1_{v_2}, \\ 
    & z_{q_3}-z^1_{q_3} + z_{q_4}-z^1_{q_4} = 1-z_{v_n}-z_{v_1}+z^1_{v_1},
    \end{split}
    \end{align}
and
\begin{align}
       & z_{s+ v^+_n} + z_{\ell(s)+ v^+_n} = z^1_{p(s)}, \; \forall s \in \E^{+}_{n-2}: s+ v^+_n, \ell(s)+ v^+_n \in \E^{+}_{n-1}\label{n20}\\
        & z_{s+ v^+_n} + z^1_{\ell(s)} = z^1_{p(s)},\; \forall s \in \E_{n-2}: s+ v^+_n \in \E^+_{n-1}, N(\ell(s))=\emptyset\label{n21}\\
       & z^1_{s} + z^1_{\ell(s)} = z^1_{p(s)}, \; \forall s \in \E^{+}_k: k \in [n-3]\setminus\{1\},{\rm or}\; k =2, N(s)=N(\ell(s))=\emptyset\label{n23}\\
    & z^1_s \geq 0,\quad \forall s \in \E_{k}: k \in [n-3], \; {\rm or} \; k=n-2, N(s)=\emptyset\label{n25}\\
    & \sum_{\substack{s \in \E^{+}_{n-2}: \\ s+v^+_n \in \E^+_{n-1}}}{z_{s+ v^+_n}}
    +\sum_{\substack{s \in \E^{+}_{n-2}:\\ N(s)=\emptyset}}{z^1_s}\leq z^1_{v_{n-1}}\label{n26}\\
    & \sum_{\substack{s \in \E^{-}_{n-2}:\\ s+v^+_n \in \E^{+}_{n-1}}}{z_{s+ v^+_n}}+\sum_{\substack{s \in \E^{-}_{n-2}: \\ N(s)=\emptyset}}{z^1_{s}} \leq z_{v_n}-z^1_{v_{n-1}}\label{n27}\\
    & \sum_{s \in \E^{+}_k}{z^1_s} \leq z^1_{v_{k+1}}, \quad \forall k \in [n-3] \setminus \{1\}\label{n28}\\
    & \sum_{s \in \E^{-}_k}{z^1_s} \leq z_{v_n}-z^1_{v_{k+1}}, \quad \forall k \in [n-3] \setminus \{1\}\label{n29}\\
    & z^1_{q_1} + z^1_{q_2} = z^1_{v_1}, \; 
    z^1_{q_1} + z^1_{q_3} = z^1_{v_2}, \; 
    z^1_{q_3} + z^1_{q_4} = z_{v_n}-z^1_{v_1},\label{n29p}
\end{align}
together with
\begin{align}
& z_{s} + z_{\ell(s)} = z_{p(s)}, \; \forall s \in \E^{+}_{n-1} \label{proj3}\\
& z_s \geq 0, \; \forall s \in \E_{n-1}\label{proj4}.
\end{align}
Note that equalities~\eqref{proj3} are present among equalities~\eqref{eqk1} and inequalities~\eqref{proj4} are present among inequalities~\eqref{eqk4}.

\medskip

In the remainder of the proof, we project out variables $z^1_v$, $v \in V \setminus \{v_n\}$, $z^1_s$, $s \in S \setminus \{s \in \E_{n-1}: N(s) \neq \emptyset\}$ in a specific order:

\paragraph{Projecting out $z^1_{v_{k}}$, $k \in \{1,\ldots,n-1\}$:}
\begin{itemize}[leftmargin=*]
    \item The variable $z^1_{v_1}$ appears only in equalities~\eqref{n19p} and~\eqref{n29p}. 
    Hence projecting out $z^1_{v_1}$ we obtain:
     \begin{align}
    &z_{q_1} + z_{q_2} = z_{v_1}\label{proj5}\\
    & z_{q_3} + z_{q_4} = 1-z_{v_1}\label{proj6}\\
    & z^1_{q_1} + z^1_{q_2}+z^1_{q_3} + z^1_{q_4} = z_{v_n}\label{int1}.
    \end{align}
    Equalities~\eqref{proj5} and~\eqref{proj6} are present among equalities~\eqref{fl}. Moreover, the above equalities
    imply:
    \begin{equation}\label{int2}
    z_{q_1}-z^1_{q_1} + z_{q_2}-z^1_{q_2}+z_{q_3}-z^1_{q_3} + z_{q_4}-z^1_{q_4} = 1-z_{v_n}
    \end{equation}
    We will use equality~\eqref{int2} to simplify our derivations.
    \item The variable $z^1_{v_2}$ appears only in equalities~\eqref{n19p} and~\eqref{n29p}.
    Hence projecting out this variable we obtain:
        \begin{align}\label{proj7}
    & z_{q_1} + z_{q_3} = z_{v_2}, 
    \end{align}
    which is present among equalities~\eqref{fl}.
    \item The variable $z^1_{v_{n-1}}$ appears only in inequalities~\eqref{n16},~\eqref{n17},~\eqref{n26}  and~\eqref{n27}.   Projecting out $z^1_{v_{n-1}}$ from~\eqref{n16} and~\eqref{n17} gives:
      \begin{equation}\label{n30}
    \sum_{s \in \E^{-}_{n-1}}{z_s}+
    \sum_{s \in \E_{n-2}: N(s) = \emptyset}{(z_s-z^1_s)}
    \leq 1-z_{v_n}, 
     \end{equation}
     while projecting out $z^1_{v_{n-1}}$ from~\eqref{n26}  and~\eqref{n27} gives:
    \begin{equation}\label{n31}
     \sum_{s \in \E^+_{n-1}}{z_s}
    +\sum_{s \in \E_{n-2}: N(s)=\emptyset}{z^1_s}\leq z_{v_{n}}.
     \end{equation}
     As we argue shortly, inequality~\eqref{n30}
     (resp.~\eqref{n31}) is implied by equality~\eqref{int2} (resp.~\eqref{int1}) and inequality~\eqref{n15} (resp.~\eqref{n25}).
     Projecting out $z^1_{v_{n-1}}$ from~\eqref{n16}
    and~\eqref{n26} gives:
    \begin{align}\label{proj8}
        & \sum_{s \in \E^{+}_{n-2}}{z_s}
    \leq z_{v_{n-1}},
    \end{align}
    which coincides with inequalities~\eqref{n1} for $k = n-2$.
    Lastly, projecting out $z^1_{v_{n-1}}$ from~\eqref{n17}
    and~\eqref{n27} we obtain:
    \begin{align}\label{proj9}
         & \sum_{s \in \E^{-}_{n-2}}{z_s}\leq 1-z_{v_{n-1}},
    \end{align}
    which coincides with inequalities~\eqref{n2} for $k = n-2$.
    \item  Variables $z^1_{v_{k+1}}$, $k \in \{2,\dots,n-3\}$, are only present in inequalities~\eqref{n18},~\eqref{n19}, ~\eqref{n28}, and~\eqref{n29}.
    Projecting out $z^1_{v_{k+1}}$, $k \in \{2,\dots,n-3\}$, from~\eqref{n18} and~\eqref{n19} gives:     
       \begin{equation*}
     \sum_{s \in \E_k}{(z_s-z^1_s)}\leq 1-z_{v_n}, \quad \forall k \in [n-3] \setminus \{1\},
     \end{equation*}
     which by equalities~\eqref{n14} can be equivalently written as:
      \begin{equation}\label{n32p}
     \sum_{s \in \E^{+}_k}{(z_{p(s)}-z^1_{p(s)})}\leq 1-z_{v_n}. \quad \forall k \in [n-3] \setminus \{1\}.
     \end{equation}
     By~\eqref{n11}--\eqref{n14}, for each $k \in [n-1] \setminus \{1\}$, we have $\sum_{s\in \E_{k+1}}{(z_s-z^1_s)} \leq \sum_{s \in \E_{k}}{(z_s-z^1_s)}$, implying that inequalities~\eqref{n30} and~\eqref{n32p} are implied by the following inequality:
     \begin{equation*}
     \sum_{s \in \E^+_2}{(z_{p(s)}-z^1_{p(s)})}\leq 1-z_{v_n},
     \end{equation*}
     which in turn is implied by inequalities~\eqref{n15} and equality~\eqref{int2}.
     Similarly, projecting out $z^1_{v_{k+1}}$, $k \in \{2,\dots,n-3\}$, from~\eqref{n28} and~\eqref{n29} gives:  
     \begin{equation*}
     \sum_{s \in \E_k}{z^1_s}\leq z_{v_n}, \quad \forall k \in [n-3] \setminus \{1\},
     \end{equation*}
     which by equalities~\eqref{n23} can be equivalently written as:
     \begin{equation}\label{n33p}
     \sum_{s \in \E^+_k}{z^1_{p(s)}}\leq z_{v_n}, \quad \forall k \in [n-3] \setminus \{1\}.
     \end{equation}
     Since for each $k \in [n-1] \setminus \{1\}$, by~\eqref{n20}--\eqref{n23}, we have $\sum_{s \in \E_{k+1}}{z^1_s} \leq \sum_{s \in \E_{k}}{z^1_s}$,  inequalities~\eqref{n31} and~\eqref{n33p} are implied by:
     \begin{equation*}
     \sum_{s \in \E^+_2}{z^1_{p(s)}}\leq z_{v_n},
     \end{equation*}
     which in turn is implied by inequalities~\eqref{n25} and equality~\eqref{int1}.
     Projecting out $z^1_{v_{k+1}}$, $k \in \{2,\dots,n-3\}$, from~\eqref{n18} and~\eqref{n28} we obtain:
     \begin{equation}\label{proj10}
     \sum_{s \in \E^+_k}{z_s} \leq z_{v_{k+1}}, \quad \forall k \in [n-3] \setminus \{1\},
     \end{equation}
     which are present among inequalities~\eqref{n1}.
     Similarly, projecting out $z^1_{v_{k+1}}$, $k \in \{2,\dots,n-3\}$, from~\eqref{n19} and~\eqref{n29}, we obtain
     \begin{equation}\label{proj11}
     \sum_{s \in \E^-_k}{z_s} \leq 1-z_{v_{k+1}}, \quad \forall k \in [n-3] \setminus \{1\},
     \end{equation}
     which are present among inequalities~\eqref{n2}.
\end{itemize}

Hence, we have shown that projecting out variables $z^1_{v_{k}}$, $k \in \{1,\ldots,n-1\}$, from~\eqref{n16}-\eqref{n19p} and~\eqref{n26}--\eqref{n29p}, we obtain inequalities~\eqref{n1}--\eqref{n2} when $k \neq n-1$, equalities~\eqref{fl}, and equalities~\eqref{int1}~\eqref{int2}.

\paragraph{Projecting out $z^1_s$, $s \in \E_{n-2}$ with $N(s) = \emptyset$:} Consider $z^1_{\bar s}$ for some $\bar s \in \E^+_{n-2}$ with $N(\bar s) = \emptyset$. 
This variable is only present in~\eqref{n12}--\eqref{n15} and~\eqref{n21}--\eqref{n25}.
Two cases arise:
\begin{itemize}[leftmargin=*]
\item If $N(\ell(\bar s)) \neq \emptyset$, then $z^1_{\bar s}$ is only present in~\eqref{n12},~\eqref{n15},~\eqref{n21},~\eqref{n25}; \ie the following system:
\begin{align*}
& z_{\bar s}-z^1_{\bar s} + z_{\ell(\bar s) + v^-_n} = z_{p(\bar s)}-z^1_{p(\bar s)}\\ 
& z_{\bar s} - z^1_{\bar s} \geq 0 \\
& z^1_{\bar s} + z_{\ell(\bar s)+v^+_n} = z^1_{p(\bar s)}\\
& z^1_{\bar s} \geq 0.
\end{align*}
Projecting out $z^1_{\bar s}$ from the above system yields:
\begin{align*}
& z_{\bar s} + z_{\ell(\bar s)} = z_{p(\bar s)}\\
& z_{\bar s} \geq 0\\
& z_{\ell(\bar s) + v^-_n} \leq  z_{p(\bar s)}-z^1_{p(\bar s)}\\ 
& z_{\ell(\bar s)+ v^+_n} \leq z^1_{p(\bar s)}  
\end{align*}
where to obtain the first equality we made use of equalities~\eqref{proj3}.
\item If $N(\ell(\bar s)) = \emptyset$, then $z^1_{\bar s}$
and $z^1_{\ell(\bar s)}$ are present only in~\eqref{n14}--\eqref{n15}, and~\eqref{n23}--\eqref{n25}; \ie the following system:
\begin{align*}
& z_{\bar s}-z^1_{\bar s} + z_{\ell(\bar s)}-z^1_{\ell(\bar s)} = z_{p(\bar s)}-z^1_{p(\bar s)},\\
& z_{\bar s} - z^1_{\bar s} \geq 0 \\
& z_{\ell(\bar s)}-z^1_{\ell(\bar s)} \geq 0\\
& z^1_{\bar s} + z^1_{\ell(\bar s)} = z^1_{p(\bar s)}, \\
& z^1_{\bar s} \geq 0 \\
& z^1_{\ell(\bar s)} \geq 0.
\end{align*}
Projecting out $z^1_{\bar s}$ and $z^1_{\ell(\bar s)}$ from the above system yields:
\begin{align*}
& z_{\bar s} + z_{\ell(\bar s)}=z_{p(\bar s)}\\
& z_{\bar s} \geq 0\\
& z_{\ell(\bar s)} \geq 0\\
& z^1_{p(\bar s)} \geq 0 \\
& z_{p(\bar s)}-z^1_{p(\bar s)} \geq 0.
\end{align*}
\end{itemize}
By a recursive application of the two steps detailed above to project out all $z^1_s$, $s \in \E_{n-2}$ with $N(s) = \emptyset$, from~\eqref{n12}--\eqref{n15} and~\eqref{n21}--\eqref{n25}, we obtain:
\begin{align*}
&z_{s + v^+_n} \leq z^1_{p(s)}, \quad \forall s \in \E_{n-2}: s + v^+_n\in \E^+_{n-1}, \; N(\ell(s)) = \emptyset, \\
&z_{s + v^-_n} \leq z_{p(s)}-z^1_{p(s)}, \quad \forall s \in \E_{n-2}: s + v^-_n\in \E^-_{n-1}, \;
N(\ell(s)) = \emptyset,\\
& z^1_{p(s)} \geq 0, \; z_{p(s)}-z^1_{p(s)} \geq 0, \quad \forall s \in \E_{n-2}: N(s)=N(\ell(s))=\emptyset
\end{align*}
together with
\begin{align}
\begin{split}\label{proj12}
& z_s + z_{\ell(s)} = z_{p(s)}, \quad \forall s \in \E^+_{n-2}\\
& z_s \geq 0, \quad \forall s \in \E_{n-2}.
\end{split}
\end{align}
Hence, to complete the projection, it suffices to project out variables $z^1_s$, $e \in \E_k$, $k \in [n-3]$, from the following system:

\begin{align}
\begin{split}\label{ssp}
&    z_{s + v^-_n} +  z_{\ell(s) + v^-_n} = z_{p(s)}-z^1_{p(s)}, \quad \forall s\in \E^+_{n-2}: s + v^-_n, \ell(s) + v^-_n \in \E^{-}_{n-1}\\
& z_{s+v^+_n} +  z_{\ell(s) + v^+_n} = z^1_{p(s)}, \quad \forall s\in \E^+_{n-2}: s + v^+_n, \ell(s) + v^+_n \in \E^{+}_{n-1}\\
&z_{s + v^-_n} \leq z_{p(s)}-z^1_{p(s)}, \quad \forall s \in \E_{n-2}: s + v^-_n\in \E^-_{n-1}, \; N(\ell(s)) = \emptyset\\
&z_{s + v^+_n} \leq z^1_{p(s)}, \quad \forall s \in \E_{n-2}: s + v^+_n\in \E^+_{n-1}, \; N(\ell(s)) = \emptyset\\
& z_s-z^1_{s} + z_{\ell(s)}-z^1_{\ell(s)} = z_{p(s)}-z^1_{p(s)}, \quad \forall s \in \E^{+}_k, k \in [n-3]\setminus\{1\}\\
& z^1_{s} + z^1_{\ell(s)} = z^1_{p(s)}, \quad \forall s \in \E^{+}_k, k \in [n-3]\setminus\{1\}\\
& z_s-z^1_s \geq 0, \quad \forall s \in \E_{k}: k \in [n-3]\\
& z^1_s \geq 0, \quad \forall  s \in \E_{k}: k \in [n-3]\\
&  z_{q_1}-z^1_{q_1} +  z_{q_2}-z^1_{q_2}+ z_{q_3}- z^1_{q_3} + z_{q_4}-z^1_{q_4}= 1-z_{v_n}\\
 &z^1_{q_1} +  z^1_{q_2}+ z^1_{q_3} + z^1_{q_4}= z_{v_n}.
\end{split}
\end{align}
For each $s \in \E^+_{n-3}$ define 
$$N^+_2(p(s)):= \{s' \in \E^+_{n-1}: p(p(s'))=s\}$$ 
and 
$$N^-_2(p(s)):= \{s' \in \E^-_{n-1}: p(p(s'))=s\}.$$ 
Note that  
$0 \leq |N^+_2(p(s))| \leq 4$ and $0 \leq |N^-_2(p(s))| \leq 4$ for all $s \in \E^+_{n-3}$.
Projecting out $z^1_{s}$, $s \in \E_{n-3}$, from system~\eqref{ssp} yields:
\begin{align*}
&z_s \geq 0, \quad \forall s \in \E_{n-3}\\
& z_s + z_{\ell(s)} = z_{p(s)}, \quad \forall s \in \E^{+}_{n-3},
\end{align*}
together with
\begin{align*}
& \sum_{s' \in N^-_2(p(s))} {z_{s'}}= z_{p(s)}-z^1_{p(s)}, \quad \forall s \in \E^{+}_{n-3}: |N^-_2(p(s))|=4\\
& \sum_{s' \in N^+_2(p(s))} {z_{s'}}= z^1_{p(s)}, \quad \forall s \in \E^{+}_{n-3}: |N^+_2(p(s))|=4\\
& \sum_{s' \in N^-_2(p(s))} {z_{s'}}\leq  z_{p(s)}-z^1_{p(s)}, \quad \forall s \in \E^{+}_{n-3}: |N^-_2(s)| < 4\\
& \sum_{s' \in N^+_2(p(e))} {z_{s'}}\leq  z^1_{p(s)}, \quad \forall s \in \E^{+}_{n-3}: |N^+_2(s)| < 4\\
& z_s-z^1_{s} + z_{\ell(s)}-z^1_{\ell(s)} = z_{p(s)}-z^1_{p(s)}, \quad \forall s \in \E^{+}_k, k \in [n-4]\setminus\{1\}\\
& z^1_{s} + z^1_{\ell(s)} = z^1_{p(s)}, \quad \forall s \in \E^{+}_k, k \in [n-4]\setminus\{1\}\\
& z_s-z^1_s \geq 0, \; \forall k \in [n-4]\\
& z^1_s \geq 0, \; \forall k \in [n-4]\\
& z^1_{q_1} +  z^1_{q_2}+ z^1_{q_3} + z^1_{q_4}= z_{v_n}\\
& z_{q_1}-z^1_{q_1} +  z_{q_2}-z^1_{q_2}+ z_{q_3}- z^1_{q_3} + z_{q_4}-z^1_{q_4}= 1-z_{v_n}
\end{align*}
By a recursive application of the above argument $n-5$ times to project out variables $z^1_s$, $s \in \E_k$, $k \in [n-4] \setminus \{1\}$, from the above system we obtain:
\begin{align}\label{sysadd1}
\begin{split}
    & \sum_{s \in \E^-_{n-1}: s \supset q_i}{z_s} \leq z_{q_i}-z^1_{q_i}, \quad \forall i\in \{1,\dots,4\}\\
    & \sum_{s \in \E^+_{n-1}: s \supset q_i}{z_s} \leq z^1_{q_i}, \quad \forall i\in \{1,\dots,4\}\\
    & z_{q_i}-z^1_{q_i} \geq 0, \quad \forall i \in \{1,\dots,4\}\\
& z^1_{q_i} \geq 0, \; \quad i \in \{1,\dots,4\}\\
& z^1_{q_1} +  z^1_{q_2}+ z^1_{q_3} + z^1_{q_4}= z_{v_n}\\
& z_{q_1}-z^1_{q_1} +  z_{q_2}-z^1_{q_2}+ z_{q_3}- z^1_{q_3} + z_{q_4}-z^1_{q_4}= 1-z_{v_n},
\end{split}
\end{align}
together with
\begin{align}\label{sysadd2}
\begin{split}
& z_s \geq 0, \quad \forall s\in \E_k, k \in [n-3]\setminus \{1\}\\
& z_s + z_{\ell(s)} = z_{p(s)}, \quad \forall e \in \E^+_k, k \in [n-3] \setminus \{1\}.
\end{split}
\end{align}
Hence it remains to project out $z^1_{q_i}$, $i \in \{1,\dots,4\}$, from system~\eqref{sysadd1}. Since $\bigcup_{i\in \{1,\dots, 4\}}{\{s \in \E^-_{n-1}: s \supset q_i\}}= \E^{-}_{n-1}$ and $\bigcup_{i\in \{1,\dots, 4\}}{\{s \in \E^+_{n-1}: s \supset q_i\}}= \E^{+}_{n-1}$, it follows that the projection of system~\eqref{sysadd1} onto the space of $z$ is given by:
\begin{align}
\begin{split}\label{final}
& \sum_{s \in \E^-_{n-1}}{z_s} \leq 1-z_{v_n}\\
& \sum_{s \in \E^+_{n-1}}{z_s} \leq z_{v_n}\\
& z_{q_i} \geq 0, \quad \forall i \in \{1,\dots,4\}\\
& z_{q_1} + z_{q_2} + z_{q_3} + z_{q_4} = 1,
\end{split}
\end{align}
where the last equality is implied by equalities~\eqref{fl}.
From~\eqref{proj3}--\eqref{proj6},~\eqref{proj7},~\eqref{proj8},~\eqref{proj9},\eqref{proj10}-\eqref{proj12},~\eqref{sysadd2}, and~\eqref{final}, it follows that the pseudo-Boolean polytope $\PBP(H)$ is defined by system~\eqref{eqk1}--\eqref{fl} and this completes the proof. 
\end{prf}

\bigskip
\noindent
\textbf{Statements and Declarations}
\bigskip

\noindent
\textbf{Competing interests:} The authors have no financial or non-financial conflicts of interests.

\noindent
\textbf{Availability of data and materials:} The manuscript has not used any data.

\bigskip
\noindent
\textbf{Funding:}
A. Del Pia is partially funded by AFOSR grant FA9550-23-1-0433. 
A. Khajavirad is in part supported by AFOSR grant FA9550-23-1-0123.
Any opinions, findings, and conclusions or recommendations expressed in this material are those of the authors and do not necessarily reflect the views of the Air Force Office of Scientific Research.

\bigskip
\noindent
\textbf{Acknowledgments:}
The authors would like to thank two anonymous referees for comments
and suggestions that improved the quality of this manuscript.

\ifthenelse {\boolean{MPA}}
{
\bibliographystyle{spmpsci}
}
{
\bibliographystyle{plainurl}
}


\end{document}